\documentclass[12pt,oneside,reqno]{amsart}
\usepackage{}
\usepackage{amssymb}
\usepackage{bbm}
\usepackage{cases}
\usepackage{amsmath}
\usepackage{graphicx}
\usepackage{mathrsfs}
\usepackage{stmaryrd}
\usepackage{amsfonts}
\usepackage{enumerate,amsmath,amssymb,amsthm}

\numberwithin{equation}{section}

\newcommand{\be}{\begin{eqnarray}}
\newcommand{\ee}{\end{eqnarray}}
\newcommand{\ce}{\begin{eqnarray*}}
\newcommand{\de}{\end{eqnarray*}}
\newtheorem{theorem}{Theorem}[section]
\newtheorem{lemma}[theorem]{Lemma}
\newtheorem{remark}[theorem]{Remark}
\newtheorem{definition}[theorem]{Definition}
\newtheorem{proposition}[theorem]{Proposition}
\newtheorem{Examples}[theorem]{Example}
\newtheorem{corollary}[theorem]{Corollary}

\renewcommand{\geq}{\geqslant}
\def\leq{\leqslant}

\newcommand{\N}{\mathbb{N}}

\newcommand{\R}{\mathbb{R}}

\newcommand{\mcl}{\mathcal}

\def\e{{\mathrm{e}}}

\def\[{{\Big[}}
\def\]{{\Big]}}
\def\<{{\langle}}
\def\>{{\rangle}}
\def\({{\Big(}}
\def\){{\Big)}}

\def\bx{{\mathbf{x}}}

\def\sgn{\mbox{\rm sgn}}

\def\dif{{\mathord{{\rm d}}}}

\def\min{{\mathord{{\rm min}}}}

\def\={&\!\!=\!\!&}
\def\bt{\begin{theorem}}
\def\et{\end{theorem}}
\def\bl{\begin{lemma}}
\def\el{\end{lemma}}
\def\br{\begin{remark}}
\def\er{\end{remark}}

\def\bd{\begin{definition}}
\def\ed{\end{definition}}
\def\bp{\begin{proposition}}
\def\ep{\end{proposition}}
\def\bc{\begin{corollary}}
\def\ec{\end{corollary}}
\def\bx{\begin{Examples}}
\def\ex{\end{Examples}}

\def\geq{\geqslant}
\def\leq{\leqslant}

\def\R{{\mathbb R}}
\allowdisplaybreaks

\begin{document}

\begin{center}
{\Large{\bf Approximation to stable law by the Lindeberg principle}}\\~\\
Peng Chen and Lihu Xu \\
{\it University of Macau} \\~\\
%\small{September 25, 2007}\\
\end{center}
{\small \noindent {\bf Abstract:}
By the Lindeberg principle, we develop in this paper an approximation to
one dimensional (possibly) asymmetric $\alpha$-stable distributions with $\alpha \in (0,2)$ in the smooth Wasserstein distance.
It is the first time that the general stable central limit theorem is proved by the Lindeberg principle, and that this theorem with $\alpha \in (0,1]$ is proved by a new method other than Fourier analysis. Our main tools are a Taylor-like expansion and a Kolmogorov forward equation.
 \\


\tableofcontents
\bigskip
\section{Motivation and main results}
Let $S_n=X_1+...+X_n$ be a sum of i.i.d. random variables whose common distribution is heavy tailed. The stable central limit theorem (CLT) tells us that under some condition \cite{Dur10}, there exists some $c_n$ such that $n^{-1/\alpha} (S_n-c_n)$ converges to an $\alpha$-stable distribution $\mu$ with $\alpha \in (0,2)$ as $n \rightarrow \infty$. There have been many works studying the convergence rate of the stable CLT in the Kolmogorov distance \cite{hall,Hall81,HauLus15,I,JuPa98,D-N}. Recently, several papers considered the convergence rate in the Wasserstein-1 distance or the smooth Wasserstein distance by Stein's method \cite{lihu,C-N-X,CGS,bluebook} or Tikhomirov-Stein's method \cite{AMPS}. However, all these works have assumed $\alpha \in (1,2)$, the more difficult regime $\alpha \in (0,1]$ is still open. This paper is along the line of the recent work by the second author \cite{lihu}, applying the Lindeberg principle to give a convergence rate for the general stable CLT with $\alpha \in (0,2)$.

Lindeberg's proof \cite{L} avoids the use of characteristic functions and gives a new and easy-to-follow way to prove the normal CLT,
it is now well known and has been well developed to study other limit theorems.
Chatterjee \cite{C-S} first applied the Lindeberg principle to identify the limiting spectral distribution of Wigner matrices with exchangeable entries, then Tao and Vu \cite{T-V} generalized this idea to prove the long standing conjecture that the universality of local eigenvalue of random matrices is determined by the first four moments of the distribution of entries. By the Lindeberg principle again, Caravenna et. al. \cite{CSZ17} obtained a general scaling limit theorem of disordered system by expanding polynomial chaos. Besides a lot of applications in random matrices \cite{BCP18,BCZ18,KnYi17,WAP17,LeSc17,Woo16,ACW16,BMP15}, the Lindeberg principle has also been applied to other research areas such as high dimensional regressions \cite{CCK14,CCK17,GNU17}, time series \cite{Gra17,WAP17,MePe16,LAP11}, statistical bootstrap \cite{Pou15,CCK17}, statistical learning \cite{KoMo11,ZCSXK13} and so on.

Let us briefly recall the application of the Lindeberg principle to the normal CLT. Let $X_{1},...,X_{n}$ be a sequence of independent random variables with mean $0$, variance $1$ and bounded third moment, to prove the normal CLT, we take a sequence of independent standard normal random variables $Y_{1},...,Y_{n}$ and compare $\frac{X_{1}+...+X_{n}}{\sqrt n}$ and $\frac{Y_{1}+...+Y_{n}}{\sqrt n}$ (note that $\frac{Y_{1}+...+Y_{n}}{\sqrt n}$ is standard normal). More precisely, denote $Z_{1}=X_{2}+...+X_{n}$, $Z_{j}=Y_{1}+...+Y_{j-1}+X_{j+1}+...+X_{n}$ for all $1<j<n$, and $Z_{n}=Y_{1}+...+Y_{n-1}$, for any bounded third order differentiable function $f$, we have
\begin{equation*}
\begin{split}
& \ \ \mathbb E\left[f\left(\frac{X_{1}+...+X_{n}}{\sqrt n}\right)-f\left(\frac{Y_{1}+...+Y_{n}}{\sqrt n}\right)\right] \\
&=\sum_{j=1}^{n}\mathbb E\left[f\left(\frac{X_{j}+Z_{j}}{\sqrt n}\right)-f\left(\frac{Z_{j}}{\sqrt n}\right)\right]-\sum_{j=1}^{n}\mathbb E\left[f\left(\frac{Y_{j}+Z_{j}}{\sqrt n}\right)-f\left(\frac{Z_{j}}{\sqrt n}\right)\right].
\end{split}
\end{equation*}
Applying third order Taylor expansion to $f$ and controlling the remainder, we finally can obtain
\begin{equation*}
\begin{split}
& \ \ \left|\mathbb E\left[f\left(\frac{X_{1}+...+X_{n}}{\sqrt n}\right)-f\left(\frac{Y_{1}+...+Y_{n}}{\sqrt n}\right)\right]\right| \le  \frac{C \sup_{x \in \mathbb R}|f'''(x)| (\mathbb E |Y|^{3}+\mathbb E |X|^{3})}{\sqrt n}.
\end{split}
\end{equation*}
The above inequality not only implies a CLT but also provides a convergence rate in a smooth Wasserstein distance.
\vskip 3mm

Although the stable CLT is one of the most important limit theorems in probability theory, surprisingly there have not been any works which apply the Lindeberg principle to prove the stable CLT. To the best of our knowledge, this paper is the first time that the general CLT is proved by the Lindeberg principle, and that this theorem with $\alpha \in (0,1]$ is proved by a new method other than Fourier analysis. Our results further provide convergence rates in the smooth Wasserstein distance, which match the best known rates in the Kolmogorov distance. Note that there is no subordination between these two distances. Although \cite{J} proved the symmetric stable CLT for $\alpha \in (0,2)$ in Mallow distance by a maximal coupling argument, its argument seems to heavily depend on the symmetry assumption and the related convergence rate is far from the best one. %To the best of our knowledge, when $\alpha \in (0,1]$, our paper is the first try which succeeds proving asymmetric stable CLT by a method other than characteristic function.

Let us give a brief discussion of our main results.  Theorem \ref{main} below provides a general convergence rates in the smooth Wasserstein distance for the stable CLT with $\alpha \in (0,2)$ when $X_{1}$ has a distribution which falls in the domain of normal attraction of stable law, while Theorem \ref{improve} further improves the rate for the case $\alpha \in (0,1]$ under a slightly stronger condition. The convergence rate in Theorem \ref{main} matches the optimal one in the Kolmogorov distance found by Hall \cite{hall}, see more details in Remark \ref{r:Rem1}.
When $X_{1}$ is out of the scope of normal attraction of stable law, we also found a convergence rate.

%When $\alpha \in (0,1)$, the rate $n^{-\min\{1,\frac 1 \alpha-1\}}$ is known to be the optimal rate in Kolmogorov distance. Some more recent literature reported a rate $n^{-1}$, which is better than the above optimal rate when $\alpha \in  (\frac 12, 1)$. In Theorem \ref{}, we clarify this contradiction.

The above direct Taylor expansion does not work for proving the stable CLT by the Lindeberg principle. Alternatively, we develop a Taylor-like expansion and use a Kolmogorov forward equation to handle the remainder. When $\alpha \in (1,2)$, this expansion is similar to that in \cite{C-N-X}. As $\alpha \in (0,1]$, due to the lack of the first moment, we need to use a truncation technique and estimate the remainder in a much more delicate way.

%Our paper first develop the Taylor-like expansion to solve the term of the sum of independent and identically distributed random variables, and use the Dynkin's formula to solve the remaineder about $\alpha-$stable law, then we prove the \emph{stable central limit theorem} entirely with the Lindeberg principle.
\ \vskip 3mm
\bd
Let $\alpha\in(0,2)$, $\sigma>0$ and $\beta\in [-1,1]$ be real numbers.
We say that $Y$ is an \emph{$\alpha$-stable} random variable with \emph{parameters $\sigma$ and $\beta$}, writing $Y\sim S_\alpha(\sigma,\beta)$, if for all $\lambda\in\R$,
    \begin{align*}
    \mathbb{E}\big[\e^{i\lambda Y}\big]=
    \left\{
\begin{array}{lll}
{\rm exp}\big\{-\sigma^\alpha|\lambda|^\alpha(1-i\,\beta\,{\rm sign}(\lambda)\,\tan\frac{\pi\alpha}2)\big\}&\quad\mbox{if $\alpha\neq 1$}\\
\\
{\rm exp}\big\{-\sigma|\lambda|(1+i\,\beta\,\frac2\pi{\rm sign}(\lambda)\,\log \lambda)\big\}&\quad\mbox{if $\alpha= 1$}
\end{array}
\right..
    \end{align*}
In particular, when $\beta=0$, we say that $Y$ has a \emph{symmetric $\alpha$-stable distribution of parameter} $\sigma$, and write $Y\sim S\alpha S(\sigma)$.
\ed
Since $Y/\sigma\sim  S_\alpha(1,\beta)$ if $Y\sim  S_\alpha(\sigma,\beta)$ (when $\alpha=1,$ we further assume $\beta=0$, which will be explained in the following), from now on we will only consider stable distributions with $\sigma=1$.

Let $(\hat{Y}_{t})_{t\geq0}$ be a one-dimensional $\alpha-$stable process, that is, a L$\acute{e}$vy process satisfying
\begin{align}  \label{e:ChfY}
\mathbb E \left[\e^{i\lambda  \hat Y_{t}}\right]\ =\ \e^{-t\psi(\lambda)},
\end{align}
where
\begin{align*}
\psi(\lambda)\ =\
    \left\{
\begin{array}{lll}
|\lambda|^\alpha(1-i\,\beta\,{\rm sign}(\lambda)\,\tan\frac{\pi\alpha}2)&\quad\mbox{if $\alpha\neq 1$}\\
\\
|\lambda|(1+i\,\beta\,\frac2\pi{\rm sign}(\lambda)\,\log \lambda)&\quad\mbox{if $\alpha= 1$}
\end{array}
\right..
\end{align*}
It is easy to see from the above $\psi$ that when $\beta \ne 0$ and $\alpha=1$, the random variable is not strictly stable because it does not have the scaling property (see e.g., \cite[Theorem 14.15]{S}). So, to rule out this singular case, we always assume \underline{$\beta=0$} when \underline{$\alpha=1$} from now on.
\vskip 3mm

We denote
$\mathcal{C}_{b}^{k}(\mathbb{R})=\{f:$ $\mathbb{R}\rightarrow\mathbb{R};$ $f,~f',~\cdots,~f^{(k)}$ are all continuous and bounded functions$\}.$  The infinitesimal generator $\mathcal{L}^{\alpha,\beta}$ of $(\tilde{Y}_{t})_{t\geq0}$ is
\bd
Let $\alpha\in(0,2)$ and $\beta\in [-1,1]$.
 For $f:\R\to\R$ in $\mcl C_{b}^2(\R)$, we define
\begin{align}\label{op}
\mathcal{L}^{\alpha,\beta}f(x)=d_\alpha \,
\int_\R \frac{f(y+x)-f(x)-y^{(\alpha)}f'(x)}{2|y|^{1+\alpha}}
\left[
(1+\beta){\bf 1}_{(0,\infty)}(y)+(1-\beta){\bf 1}_{(-\infty,0)}(y)
\right]
 dy,
\end{align}
where $d_\alpha = \left(
\int_0^\infty \frac{1-\cos y}{y^{1+\alpha}}dy
\right)^{-1}$ and
\begin{align*}
y^{(\alpha)}=
\begin{cases}
y, &\alpha\in(1,2),\\
y{\bf 1}_{(-1,1)}(y), &\alpha=1,\\
0, &\alpha\in(0,1).
\end{cases}
\end{align*}
When $\beta=0,$ $\mathcal{L}^{\alpha,0}$ reduces to the fractional Laplacian $\Delta^\frac\alpha2$ of order $\alpha/2$.
\ed
Recall that the Wasserstein-1 distance between two probability measures $\mu_{1}$ and $\mu_{2}$ is defined by
$$
W_1(\mu_{1},\mu_{2})\ =\ \sup_{h \in {\rm Lip}(1)} |\mu_{1}(h)-\mu_{2}(h)|,
$$
where ${\rm Lip}(1)=\{h: \R \rightarrow \R; \ |h(y)-h(x)| \le |y-x|\}$ and $\mu_i(h)=\int h(x) \mu_i(\dif x)$ for $i=1,2$. The Kolmogorov distance is
$$d_{{\rm Kol}}(\mu_{1}, \mu_{2})=\sup_{x \in \R} |\mu_1((-\infty,x])-\mu_2((-\infty,x])|.$$
The smooth Wasserstein distance of order $k \in \N$ (\cite[(4.1)]{AH}) is
$$d_{\mcl W_k} (\mu_1, \mu_2) \ = \ \sup_{h \in \mcl H_{k}} \left|\mu_1 (h)-\mu_2(h)\right|,$$
$\mcl H_{k}$ is the set of all bounded $k$-th order differentiable functions $h$ such that $\|h^{{(j)}}\| \le 1$ for $j=0,1,...,k$.
\ \ \\

Before giving the first main theorem, we first recall the definition of normal attraction of a stable law of exponent $\alpha$.

\bd
If $X$ has a distribution function of the form
\begin{equation}\label{px}
F_X(x)=\big(1-\frac{A+\epsilon(x)}{x^{\alpha}}(1+\beta)\big){\bf 1}_{[0,\infty)}(x)+\frac{A+\epsilon(x)}{|x|^{\alpha}}(1-\beta){\bf 1}_{(-\infty,0)}(x)\big),
\end{equation}
where $\alpha\in(0,2)$, $A>0$, $\beta\in[-1,1]$ and  $\epsilon: \R\to\R$ is a bounded measurable function vanishing at $\pm\infty$, then we say that $X$ is in the domain $\mathcal{D}_\alpha$ of normal attraction of a stable law of exponent $\alpha$.
\ed

Let $X_1,X_2,\ldots$ be i.i.d. random variables with common distribution function $F_X$ defined by \eqref{px}. Set $\sigma=\left(A\alpha\int_{\mathbb{R}}\frac{1-\cos y}{|y|^{1+\alpha}}\dif y\right)^\frac1\alpha$ and
\begin{align}\label{sum}
S_{n}=\frac{1}{\sigma}n^{-\frac{1}{\alpha}}
\begin{cases}
X_{1}+\cdots+X_{n}-n\mathbb{E}[X_{1}], \quad &\alpha\in(1,2),\\
X_{1}+\cdots+X_{n}-n\mathbb{E}\big[X_{1}{\bf 1}_{(0,\sigma n^{\frac{1}{\alpha}})}(|X_{1}|)\big], &\alpha=1,\\
X_{1}+\cdots+X_{n}, &\alpha\in(0,1).
\end{cases}
\end{align}
We have the following theorem:
%According to (\ref{px}) and (\ref{bound}), we can give the following main results.
%Let us have a discussion on the case $\alpha=1$.

\bt\label{main}
Let $Y$ have the distribution $S_{\alpha}(1,\beta)$ with $\alpha \in (0,2)$, $\beta \in [-1,1]$, and $\beta=0$ in the case $\alpha=1$. Assume that
there exist some constants $K>0$ and $\gamma\geq0$ such that
\begin{align}\label{bound}
|\epsilon(x)|\leq\frac{K}{|x|^{\gamma}},\quad x\neq0.
\end{align}
Then for any $f\in\mathcal{C}_{b}^{3}(\mathbb{R}),$ there exists $c_{\alpha,A,K}$ (that can be made explicit) depending only on $\alpha$, $A$ and $K$ such that

\noindent i) When $\alpha\in(1,2),$ we have
\begin{align*}
\big|\mathbb{E}[f(S_{n})]-\mathbb{E}[f(Y)]\big|\leq c_{\alpha,A}&\big(\|f'\|_{\infty}+\|f''\|_{\infty}+\|f^{(3)}\|_{\infty}\big)\\
&\cdot\Big[
n^{\frac{\alpha-2}{\alpha}}\Big(1+\int_{-\sigma\,n^{\frac1\alpha}}^{\sigma\,n^{\frac1\alpha}}
\frac{|\epsilon(x)|}{|x|^{\alpha-1}}dx\Big)
+\sup_{|x|\geq \sigma\,n^{\frac1\alpha}}|\epsilon(x)|
\Big].
\end{align*}
ii) When $\alpha=1,$ $\beta=0$ and $\gamma>0,$ we have
\begin{align*}
\big|\mathbb{E}[f(S_{n})]-\mathbb{E}[f(Y)]\big|&\leq c_{\alpha,A}\big(\|f\|_{\infty}+\|f'\|_{\infty}+\|f''\|_{\infty}+\|f^{(3)}\|_{\infty}\big)\\
&\cdot\Big(n^{-1}\log n+\sup_{|x|\geq \sigma n}|\epsilon(x)|+n^{-1}\int_{-\sigma n}^{\sigma n}|\epsilon(x)|dx+\int_{\sigma n}^{\infty}\frac{|\epsilon(x)|+|\epsilon(-x)|}{x}dx\Big).
\end{align*}
iii) When $\alpha\in(0,1),$ we have
\begin{align*}
\big|\mathbb{E}[f(S_{n})]-\mathbb{E}[f(Y)]\big|\leq\|f'\|_{\infty}\mathcal{R}_{\alpha,\beta}(n)&+c_{\alpha,A,K}\big(\|f\|_{\infty}+\|f'\|_{\infty}+\|f''\|_{\infty}\big)\\ &\cdot\Big(n^{-1}+n^{\frac{\alpha-2}{\alpha}}\int_{-\sigma n^{\frac{1}{\alpha}}}^{\sigma n^{\frac{1}{\alpha}}}\frac{|\epsilon(x)|}{|x|^{\alpha-1}}dx+\mathcal{R}_{\alpha,\beta,\gamma}(n)
\Big),
\end{align*}
where
\begin{align}\label{rm}
\mathcal{R}_{\alpha,\beta}(n)=\frac{n^{\frac{\alpha-1}{\alpha}}}{\sigma}\Big|(1+\beta)\int_{0}^{\sigma n^{\frac{1}{\alpha}}}\frac{\epsilon(x)}{x^{\alpha}}\dif x-(1-\beta)\int_{0}^{\sigma n^{\frac{1}{\alpha}}}\frac{\epsilon(-x)}{x^{\alpha}}\dif x\Big|
\end{align}
and
\begin{align*}
\mathcal{R}_{\alpha,\beta,\gamma}(n)=
\begin{cases}
\sup_{|x|\geq\sigma n^{\frac{1}{\alpha}}}|\epsilon(x)|+n^{\frac{\alpha-1}{\alpha}}\int_{\sigma n^{\frac{1}{\alpha}}}^{\infty}\frac{|\epsilon(x)|+|\epsilon(-x)|}{x^{\alpha}}dx, &\gamma\in(1-\alpha,\infty),\\
\sup_{|x|\geq\sigma n^{\frac{1}{\alpha}}}|\epsilon(x)|^{\alpha}, &\gamma\in[0,1-\alpha].
\end{cases}
\end{align*}
\et
\ \

\begin{remark}  \label{r:Rem1}
From Theorem \ref{main} and the definition of $d_{\mcl W_k}$, we have
$$d_{\mcl W_3}(S_n,Y) \le \ c_{\alpha,A}\Big[
n^{\frac{\alpha-2}{\alpha}}\Big(1+\int_{-\sigma\,n^{\frac1\alpha}}^{\sigma\,n^{\frac1\alpha}}
\frac{|\epsilon(x)|}{|x|^{\alpha-1}}dx\Big)
+\sup_{|x|\geq \sigma\,n^{\frac1\alpha}}|\epsilon(x)|
\Big], \ \ \ \ {\rm for} \ \ \alpha \in (1,2); $$
\begin{align*}
d_{\mcl W_3}(S_n,Y) \le \ c_{\alpha,A}
\Big[n^{-1}\log n+\sup_{|x|\geq \sigma n}|\epsilon(x)|+n^{-1}\int_{-\sigma n}^{\sigma n}|\epsilon(x)|dx+\int_{\sigma n}^{\infty}\frac{|\epsilon(x)|+|\epsilon(-x)|}{x}dx\Big],
\end{align*}
for $\alpha=1,$ $\beta=0$ and $\gamma\in(0,\infty);$
\begin{align*}
d_{\mcl W_2}(S_n,Y) \le \ c_{\alpha,A,K}\Big[n^{-1}+n^{\frac{\alpha-2}{\alpha}}\int_{-\sigma n^{\frac{1}{\alpha}}}^{\sigma n^{\frac{1}{\alpha}}}\frac{|\epsilon(x)|}{|x|^{\alpha-1}}dx+\mathcal{R}_{\alpha,\beta,\gamma}(n)+\mathcal{R}_{\alpha,\beta}(n)
\Big], \ \ \ \ {\rm for} \ \ \alpha \in (0,1).
\end{align*}

Let us compare our rates with those in \cite{hall} and \cite{KuKe00}, where \cite{KuKe00} only studied the convergence rate in the symmetric distribution case. We consider the case that $\gamma>0$ is sufficiently large so that $\epsilon(x)$ is negligible. This case corresponds to that $D(x)$ and $S(x)$ in \cite{hall} are both negligible. When $\alpha \in (1,2)$, our rate is $n^{\frac{\alpha-2}{\alpha}}$, the same as those in \cite{hall} and \cite{KuKe00}. When $\alpha=1$, our rate is $n^{-1} \log n$, consistent with that in \cite{KuKe00}, but \cite{hall} gives a worse $n^{-1} (\log n)^2$. When $\alpha \in (0,1)$, if $\beta \ne 0$ we have $\mathcal R(\alpha,\beta)=O(n^{\frac{\alpha-1}{\alpha}})$ and get a rate $n^{-1} \vee n^{\frac{\alpha-1}{\alpha}}$, it is the same as that in \cite{hall}. As $\beta=0$, we have $\mathcal R(\alpha,\beta)=0$ and thus the rate is $n^{-1}$, consistent with the result in \cite{KuKe00}. However, the rate in \cite{hall} does not change with $\beta$. Hence, our estimate is more sensitive with respect to $\beta$.
\end{remark}

If $\epsilon(x)\rightarrow0$ as $x\rightarrow\pm\infty$, then we have $\big|\mathbb{E}[f(S_{n})]-\mathbb{E}[f(Y)]\big| \rightarrow 0$ from the previous theorem. By the same argument as the proof of \cite[Corollary I.1]{C-D}, we get the following upper bound for the Kolmogorov distance, see the proof in Appendix.
\begin{corollary}  \label{c:SCLT}
Keep the same notations and assumptions as in Theorem \ref{main}. Then we have
\begin{align*}
d_{{\rm Kol}}(S_n,Y):=\sup_{x\in\mathbb{R}}\Big|\mathbb{P}(S_{n}\leq x)-P(Y\leq x)\Big|=
\begin{cases}
O\Big(\big(d_{\mcl W_3}(S_n,Y)\big)^{\frac{1}{4}}\Big), &\alpha\in[1,2),\\
O\Big(\big(d_{\mcl W_2}(S_n,Y)\big)^{\frac{1}{3}}\Big), &\alpha\in(0,1).
\end{cases}
\end{align*}
\end{corollary}

Our next result gives an improved upper bound on $\big|\mathbb{E}[f(S_{n})]-\mathbb{E}[f(Y)]\big|$ for $\alpha\in(0,1],$ under slightly more restrictive conditions (see, e.g., \cite[Theorem 2]{hall}).
\bt\label{improve}
Consider $\alpha\in(0,1].$ Keep the same notations and assumptions as in Theorem \ref{main}. In addition, we further assume $\frac{\epsilon(x)}{x^{\alpha}}{\bf 1}_{(0,\infty)}(x)$ and $\frac{\epsilon(x)}{|x|^{\alpha}}{\bf 1}_{(-\infty,0)}(x)$ are ultimately monotone (that is, there exists $x_{0}>0$ such that $\frac{\epsilon(x)}{x^{\alpha}}{\bf 1}_{(0,\infty)}(x)$ and $\frac{\epsilon(x)}{|x|^{\alpha}}{\bf 1}_{(-\infty,0)}(x)$ are monotone for any $|x|>x_{0}$). Then there exists $\hat{c}_{\alpha,A}$ (that can be made explicit) depending only on $\alpha$ and $A$ such that \\
i) When $\alpha=1$ and $\beta=0,$ we have
\begin{align*}
\big|\mathbb{E}[f(S_{n})]-\mathbb{E}[f(Y)]\big|\leq  \hat{c}_{\alpha,A}&\big(\|f\|_{\infty}+\|f'\|_{\infty}+\|f''\|_{\infty}+\|f^{(3)}\|_{\infty}\big)\\
&\cdot\big(n^{-1}\log n+n^{-1}\int_{-\sigma n}^{\sigma n}|\epsilon(x)|\dif x+\sup_{|x|\geq\sigma n}|\epsilon(x)|\big).
\end{align*}
ii) When $\alpha\in(0,1),$ we have
\begin{align*}
\big|\mathbb{E}[f(S_{n})]-\mathbb{E}[f(Y)]\big|\leq\|f'\|_{\infty}\mathcal{R}_{\alpha,\beta}(n)&+\hat{c}_{\alpha,A}\big(\|f\|_{\infty}+\|f'\|_{\infty}+\|f''\|_{\infty}\big)\\
&\cdot\big(n^{-1}+n^{\frac{\alpha-2}{\alpha}}\big(1+\int_{-\sigma n^{\frac{1}{\alpha}}}^{\sigma n^{\frac{1}{\alpha}}}\frac{|\epsilon(x)|}{|x|^{\alpha-1}}\dif x\big)+\sup_{|x|\geq\sigma n^{\frac{1}{\alpha}}}|\epsilon(x)|\big),
\end{align*}
where $\mathcal{R}_{\alpha,\beta}(n)$ is defined by (\ref{rm}).
\et

\medskip
Under the condition that the function $\epsilon$ is a bounded measurable function vanishing at $\pm\infty,$ we obtain the stable CLT in Corollary \ref{c:SCLT}, but this condition is not necessary for the stable CLT. By slightly modifying the approach leading to Theorem \ref{main}, we can also consider the case where $\epsilon$ is a slowly varying function diverging at infinity. Because it would be too technical to state such result at a great level of generality, we prefer to illustrate an explicit example in \cite{lihu}, for which our methodology still allows to conclude.

\medskip

{\it \underline{Example}: slowly varying tails in \cite{lihu}}. We consider
\begin{align*}
P(X>x)=K_{0}\frac{\big(\log|x|\big)^{\delta}}{|x|^{\alpha}},\quad x\geq e,\qquad P(X<x)=K_{0}\frac{\big(\log|x|\big)^{\delta}}{|x|^{\alpha}},\quad x\leq-e,
\end{align*}
where $K_{0}>0,$ $\alpha\in(0,2)$ and $\delta\in\mathbb{R}.$ Define $\gamma_{n}=\inf\{x>0:\mathbb{P}(|X|>x)\leq\frac{1}{n}\}$ and $\tilde{\sigma}=(\frac{\alpha}{d_{\alpha}})^{\frac{1}{\alpha}}.$ Let $\tilde{S}_{n}=\frac{1}{\tilde{\sigma}\gamma_{n}}\big(X_{1}+\cdots+X_{n}\big),$ and consider $\tilde{Y}\sim S\alpha S(1)$.
We can deduce from a suitable modification of Theorem \ref{main} (see Section 4)  that for any $f\in\mathcal{C}_{b}^{3}(\mathbb{R}),$
\begin{align}\label{ex}
\big|\mathbb{E}[f(\tilde{S}_{n})]-\mathbb{E}[f(\tilde{Y})]\big|=
\begin{cases}
O\big((\log n)^{-1+\frac{1}{\alpha+1}}\big), &\alpha\in(0,1),\\
O\big((\log n)^{\theta-\frac{1}{2}}\big), &\alpha=1,\\
O\big((\log n)^{-1+\frac{1}{\alpha}}\big), &\alpha\in(1,2),
\end{cases}
\end{align}
for some very small $\theta>0.$ In particular, when $\delta\in(0,1],$ we have
\begin{align}\label{ex1}
\big|\mathbb{E}[f(\tilde{S}_{n})]-\mathbb{E}[f(\tilde{Y})]\big|=
\begin{cases}
O\big((\log n)^{-1+\frac{1}{\alpha+1}}\wedge(\log n)^{-\delta}\big), &\alpha\in(0,1),\\
O\big((\log n)^{\theta-\frac{1}{2}}\wedge(\log n)^{-\delta}\big), &\alpha=1,\\
O\big((\log n)^{-1+\frac{1}{\alpha}}\wedge(\log n)^{-\delta}\big), &\alpha\in(1,2),
\end{cases}
\end{align}
for some very small $\theta>0,$ here $a\wedge b=\min\{a,b\}.$

\br
For the case $\alpha\in(1,2),$ our convergence rate is $(\log n)^{-1+\frac{1}{\alpha}},$ which is consistent with the result in \cite{C-N-X}.
\er
\medskip

The rest of the paper is organized as follows. In Section \ref{idea}, we give a short proof of the Theorem \ref{main} in a special case to illustrate the main idea. In Section \ref{useful properties}, we first give some useful properties of the operator $\mathcal{L}^{\alpha,\beta}$ and asymmetric $\alpha-$stable process, then we develop the Taylor-like expansion. In Section \ref{proof}, we extend the Lindeberg principle to the asymmetric $\alpha-$stable distributions and provide the proofs of Theorem \ref{main} and Theorem \ref{improve}. In Section \ref{difficult}, we will focus on the proofs of (\ref{ex}) and (\ref{ex1}). We will prove Corollary \ref{c:SCLT} in the last section.

Let $(X_t)_{t \ge 0}$ be a time homogeneous $\R$-valued stochastic process with the infinitesimal generator $\mathcal G$. For a bounded continuous function $f$, define $P_t f(x)=\mathbb{E} [f(X_t)]$ with $X_0=x \in \R$, $(P_{t})_{t\geq0}$ is a semigroup on the space of bounded continuous functions. The Kolmogorov forward equation \cite[Section 4.6 and (A.9.4)]{B-K} tells us that as long as $f$ is in the domain of $\mathcal G$, we have
$$P_t f(x)\ =\ f(x)+\int_0^t P_s \mathcal G f(x) \dif s,$$
which reads in the language of stochastic process as
\begin{align}  \label{e:SpDyk}
\mathbb{E}_{X}[f(X_{t})]-f(X)=\int_0^t \mathbb{E}_{X}\Big[\mathcal{G}f(X_{s})\Big]\dif s,
\end{align}
here $X_{0}=X$ and we assume that all the functions $f$ in the sequel are in the domain of $\mathcal G$.

\section{A short proof of Theorem \ref{main} in a special case to illustrate the main idea}\label{idea}
As mentioned in the introduction, we shall use a Taylor-like expansion and (\ref{e:SpDyk}) to prove Theorem \ref{main}. Before we go into the details, let us give its short proof in the symmetric Pareto distribution case (see, e.g., \cite{D-N}) to illustrate how these two tools will work.

Assume that $X_1,X_2,\ldots$ are independent copies drawn from the \emph{Pareto law of index $\alpha\in(0,2)$}, that is, suppose that the
common density is
$$
p(x)=\frac{\alpha}{2}|x|^{-(1+\alpha)}{\bf 1}_{[1,\infty)}(|x|).
$$
Consider $\sigma=\left(\frac{\alpha}2\int_\R \frac{1-\cos y}{|y|^{1+\alpha}}dy\right)^{\frac1\alpha}$ and $S_n=\frac{1}\sigma n^{-\frac1\alpha}(X_1+\ldots+X_n)$. Let $Y$ be an $\alpha$-stable random variable with characteristic function $e^{-|\lambda|^\alpha}$. Set
$$
Z_{i}=Y_{1}+\cdots+Y_{i-1}+\frac{1}{\sigma}X_{i+1}+\cdots+\frac{1}{\sigma}X_{n},\quad 1\leq i\leq n,
$$
where $Y_{1},Y_{2},\cdots$ are independent copies of $Y$, we know $\frac{Y_{1}+\cdots+Y_{n}}{n^{\frac{1}{\alpha}}}$ has the same distribution as $Y$. Then,
\begin{align*}
\mathbb{E}\big[f(S_{n})\big]-\mathbb{E}\big[f(Y)\big]%&=\sum_{i=1}^{n}\mathbb{E}\Big[f\Big(\frac{X_{1}+\cdots+X_{n}}{\sigma n^{\frac{1}{\alpha}}}\Big)-f\Big(\frac{Y_{1}+\cdots+Y_{n}}{n^{\frac{1}{\alpha}}}\Big)\Big]\\
&=\sum_{i=1}^{n}\mathbb{E}\Big[f\Big(\frac{Z_{i}}{n^{\frac{1}{\alpha}}}+\frac{X_{i}}{\sigma n^{\frac{1}{\alpha}}}\Big)-f\Big(\frac{Z_{i}}{n^{\frac{1}{\alpha}}}\Big)\Big]-\sum_{i=1}^{n}\mathbb{E}\Big[f\Big(\frac{Z_{i}}{n^{\frac{1}{\alpha}}}+\frac{Y_{i}}{n^{\frac{1}{\alpha}}}\Big)
-f\Big(\frac{Z_{i}}{n^{\frac{1}{\alpha}}}\Big)\Big].
\end{align*}
For the first term, since $Z_{i}$ and $X_{i}$ are independent, we have
\begin{align*}
\mathbb{E}\Big[f\Big(\frac{Z_{i}}{n^{\frac{1}{\alpha}}}+\frac{X_{i}}{\sigma n^{\frac{1}{\alpha}}}\Big)-f\Big(\frac{Z_{i}}{n^{\frac{1}{\alpha}}}\Big)\Big]&=\mathbb{E}\Big[\frac{\alpha}{2}\int_{|x|\geq1}\frac{f\Big(\frac{Z_{i}}{n^{\frac{1}{\alpha}}}+\frac{x}{\sigma n^{\frac{1}{\alpha}}}\Big)-f\Big(\frac{Z_{i}}{n^{\frac{1}{\alpha}}}\Big)}{|x|^{\alpha+1}}\dif x\Big]\\
&=n^{-1}\mathbb{E}\Big[\frac{d_{\alpha}}{2}\int_{|y|\geq\sigma^{-1}n^{-\frac{1}{\alpha}}}\frac{f\Big(\frac{Z_{i}}{n^{\frac{1}{\alpha}}}+y\Big)-f\Big(\frac{Z_{i}}{n^{\frac{1}{\alpha}}}\Big)}{|y|^{\alpha+1}}\dif y\Big]\\
&=n^{-1}\mathbb{E}\Big[\Delta^{\frac{\alpha}{2}}f\Big(\frac{Z_{i}}{n^{\frac{1}{\alpha}}}\Big)-\frac{d_{\alpha}}{2}\int_{-\sigma^{-1}n^{-\frac{1}{\alpha}}}^{\sigma^{-1}n^{-\frac{1}{\alpha}}}\frac{f\Big(\frac{Z_{i}}{n^{\frac{1}{\alpha}}}+y\Big)-f\Big(\frac{Z_{i}}{n^{\frac{1}{\alpha}}}\Big)}{|y|^{\alpha+1}}\dif y\Big],
\end{align*}
and
\begin{align*}
\quad\int_{-\sigma^{-1}n^{-\frac{1}{\alpha}}}^{\sigma^{-1}n^{-\frac{1}{\alpha}}}\frac{f\Big(\frac{Z_{i}}{n^{\frac{1}{\alpha}}}+y\Big)
-f\Big(\frac{Z_{i}}{n^{\frac{1}{\alpha}}}\Big)}{|y|^{\alpha+1}}\dif y&=\int_{-\sigma^{-1}n^{-\frac{1}{\alpha}}}^{\sigma^{-1}n^{-\frac{1}{\alpha}}}\frac{f\Big(\frac{Z_{i}}{n^{\frac{1}{\alpha}}}+y\Big)
-f\Big(\frac{Z_{i}}{n^{\frac{1}{\alpha}}}\Big)-yf'\Big(\frac{Z_{i}}{n^{\frac{1}{\alpha}}}\Big)}{|y|^{\alpha+1}}\dif y\\
&=\int_{-\sigma^{-1}n^{-\frac{1}{\alpha}}}^{\sigma^{-1}n^{-\frac{1}{\alpha}}}\int_{0}^{1}\frac{y\Big[f'\Big(\frac{Z_{i}}{n^{\frac{1}{\alpha}}}+yt\Big)
-f'\Big(\frac{Z_{i}}{n^{\frac{1}{\alpha}}}\Big)\Big]}{|y|^{\alpha+1}}\dif t\dif y \\
&=\int_{-\sigma^{-1}n^{-\frac{1}{\alpha}}}^{\sigma^{-1}n^{-\frac{1}{\alpha}}}\int_{0}^{1}\int_{0}^{1}\frac{ty^{2}f''\Big(\frac{Z_{i}}{n^{\frac{1}{\alpha}}}+yt\theta\Big)}{|y|^{\alpha+1}}\dif\theta\dif t\dif y,
\end{align*}
which further gives
\begin{align*}
\left|n^{-1}\mathbb{E}\Big[\frac{d_{\alpha}}{2}\int_{|y|<\sigma^{-1}n^{-\frac{1}{\alpha}}}\frac{f\Big(\frac{Z_{i}}{n^{\frac{1}{\alpha}}}+y\Big)
-f\Big(\frac{Z_{i}}{n^{\frac{1}{\alpha}}}\Big)}{|y|^{\alpha+1}}\dif y\Big]\right| \le  n^{-1}\frac{d_{\alpha}\|f''\|_{\infty}}{\sigma^{2-\alpha}(2-\alpha)}n^{\frac{\alpha-2}{\alpha}}.
\end{align*}
For the second term, since $Y_{i}$ has the same distribution as $\hat Y_1$, where $(\hat Y_t)_{t \ge 0}$ is a symmetric $\alpha$-stable process with characteristic function $e^{-t |\lambda|^\alpha}$. We have by (\ref{e:SpDyk})
 \begin{align*}
 \mathbb{E}\Big[f\Big(\frac{Z_{i}}{n^{\frac{1}{\alpha}}}+\frac{Y_{i}}{n^{\frac{1}{\alpha}}}\Big)
-f\Big(\frac{Z_{i}}{n^{\frac{1}{\alpha}}}\Big)\Big]=\int_{0}^{1}\mathbb{E}\Big[\Delta^{\frac{\alpha}{2}}_{\hat{Y}_{s}}f\Big(\frac{Z_{i}}{n^{\frac{1}{\alpha}}}
+\frac{\hat{Y}_{s}}{n^{\frac{1}{\alpha}}}\Big)\Big] \dif s =n^{-1}\int_{0}^{1}\mathbb{E}\Big[\Delta^{\frac{\alpha}{2}}f\Big(\frac{Z_{i}}{n^{\frac{1}{\alpha}}}
+\frac{\hat{Y}_{s}}{n^{\frac{1}{\alpha}}}\Big)\Big] \dif s.
 \end{align*}
 Hence,
 \begin{align*}
 \left|\mathbb{E}\big[f(S_{n})\big]-\mathbb{E}\big[f(Y)\big] \right|  \le \frac 1n \sum_{i=1}^n\left|\int_0^1  \mathbb{E}\Big[\Delta^{\frac{\alpha}{2}}f\Big(\frac{Z_{i}}{n^{\frac{1}{\alpha}}}
 +\frac{\hat{Y}_{s}}{n^{\frac{1}{\alpha}}}\Big)-\Delta^{\frac{\alpha}{2}}f\Big(\frac{Z_{i}}{n^{\frac{1}{\alpha}}}\Big)\Big] \dif s\right|+\frac{d_{\alpha}\|f''\|_{\infty}}{\sigma^{2-\alpha}(2-\alpha)}n^{\frac{\alpha-2}{\alpha}}.
 \end{align*}
Applying \eqref{22} for $\alpha \in (1,2)$ and \eqref{11} for $\alpha \in (0,1]$ with a straightforward calculation, we immediately obtain
\begin{align*}
\big|\mathbb{E}\big[f(S_{n})\big]-\mathbb{E}\big[f(Y)\big]\big|\leq c_{\alpha}
\begin{cases}
\big(\|f'\|_{\infty}+\|f''\|_{\infty}+\|f'''\|_{\infty}\big)n^{\frac{\alpha-2}{\alpha}}, & \alpha \in (1,2),\\
\big(\|f\|_{\infty}+\|f'\|_{\infty}+\|f''\|_{\infty}+\|f'''\|_{\infty}\big)n^{-1}\log n,\quad &\alpha=1,\\
\big(\|f\|_{\infty}+\|f'\|_{\infty}+\|f''\|_{\infty}\big)n^{-1},\quad &\alpha\in(0,1).
\end{cases}
 \end{align*}
 Note that the proofs of \eqref{22} and \eqref{11} are much simpler in the symmetric distribution case.
\section{Preliminaries of stable processes and nonlocal operators}\label{useful properties}
Let us first give the following heat kernel estimates, which will be used in the analysis in this section.
\bl\label{g1}
Let $(\hat{Y}_{t})_{t\geq0}$ be defined by \eqref{e:ChfY}. Then the distribution of $\hat{Y}_{t}$ has a density $p(t,x)$ for all $t>0$. Moreover, for any $t\in(0,1),$ there exists a constant $C_{\alpha}>1$ such that
\begin{align}\label{G1}
p(t,x)\leq C_{\alpha}t^{-\frac{1}{\alpha}}\Big(1\wedge\frac{t^{\frac{\alpha+1}{\alpha}}}{|x|^{\alpha+1}}\Big).
\end{align}
\el
\begin{proof}
When $\beta\in(-1,1),$ (\ref{G1}) follows from \cite[Theorem 1.1 (iii)]{C-Z}. When $\beta=1$ or $-1,$ we first have by scaling property,
 \begin{align}\label{G1}
p(t,x)=t^{-\frac{1}{\alpha}}p(1,t^{-\frac{1}{\alpha}}x),
\end{align}
here the $p(1,x)$ is also the density of one-dimensional stable distribution. Hence, according to Proposition 2.2, (2.4.8) and (2.5.4) in \cite[Chapter 2]{V-Z}, we have
\begin{align}\label{density1}
\lim_{|x|\rightarrow\infty}p(1,x)|x|^{1+\alpha}=C_{\alpha},
\end{align}
for some constant $C_{\alpha}$ depends on $\alpha.$ In addition, by the inverse of Fourier transform, we have
\begin{align}\label{density2}
p(x)=\frac{1}{2\pi}\int_{\mathbb{R}}\e^{-\psi(\lambda)}\e^{-ix\lambda}\dif\lambda&=\frac{1}{2\pi}\int_{\mathbb{R}}\e^{-|\lambda|^{\alpha}}\cos[|\lambda|^{\alpha}\beta\tan\frac{\pi\alpha}{2}\sgn(\lambda)
-x\lambda]\dif\lambda\leq\frac{\Gamma(\frac{1}{\alpha})}{\pi\alpha},
\end{align}
here $\Gamma(\cdot)$ is the gamma function. Therefore, combining (\ref{density1}) and (\ref{density2}), we obtain the desired result.
\end{proof}

\subsection{Estimates for the operator $\mathcal{L}^{\alpha,\beta}$}
By the definition of operator $\mathcal{L}^{\alpha,\beta},$ we can first get the following Proposition.
\bp\label{p}
For any $f\in\mathcal{C}_{b}^{3}(\mathbb{R}),$ and $x,z\in\mathbb{R},$ we have
\begin{align*}
|\mathcal{L}^{\alpha,\beta}f(x+z)-\mathcal{L}^{\alpha,\beta}f(x)|\leq D_{\alpha}|z|,
\end{align*}
where
\begin{align*}
D_{\alpha}=
\begin{cases}
\frac{2d_{\alpha}\|f''\|_{\infty}}{\alpha-1}+\frac{d_{\alpha}\|f^{(3)}\|_{\infty}}{2(2-\alpha)},\quad&\alpha\in(1,2),\\
2d_{\alpha}\|f'\|_{\infty}+\frac{d_{\alpha}}{2}\|f^{(3)}\|_{\infty},\qquad &\alpha=1,\\
\frac{2}{\alpha}d_{\alpha}\|f'\|_{\infty}+\frac{1}{1-\alpha}d_{\alpha}\|f''\|_{\infty}, \quad&\alpha\in(0,1).
\end{cases}
\end{align*}
\ep
\begin{proof}
For convenience, we denote
$$
I_{\beta}(y)=(1+\beta){\bf 1}_{(0,1]}(y)+(1-\beta){\bf 1}_{[-1,0)}(y),\quad I^{\beta}(y)=(1+\beta){\bf 1}_{(1,\infty)}(y)+(1-\beta){\bf 1}_{(-\infty,-1)}(y).
$$
1. When $\alpha\in(1,2),$ we have by (\ref{op})
\begin{align*}
\frac{1}{d_{\alpha}}\mathcal{L}^{\alpha,\beta}f(x)&=\int_{-\infty}^{\infty}\frac{f(y+x)-f(x)-yf'(x)}{2|y|^{1+\alpha}}I^{\beta}(y)\dif y+\int_{-\infty}^{\infty}\frac{f(y+x)-f(x)-yf'(x)}{2|y|^{1+\alpha}}I_{\beta}(y)\dif y\\
&=\int_{-\infty}^{\infty}\int_{0}^{1}\frac{yf'(x+ty)-yf'(x)}{2|y|^{1+\alpha}}
I^{\beta}(y)\dif t\dif y+\int_{-\infty}^{\infty}\int_{0}^{1}\int_{0}^{1}\frac{tf''(x+uty)}{2|y|^{\alpha-1}}
I_{\beta}(y)\dif u\dif t\dif y,
\end{align*}
it follows that
\begin{align*}
&\quad\frac{1}{d_{\alpha}}\big|\mathcal{L}^{\alpha,\beta}f(x+z)-\mathcal{L}^{\alpha,\beta}f(x)\big|\\
&\leq \Big|\int_{-\infty}^{\infty}\int_{0}^{1}\frac{y[f'(x+z+ty)-f'(x+z)-f'(x+ty)+f'(x)]}{2|y|^{1+\alpha}}
I^{\beta}(y)\dif t\dif y\Big|\\
&\quad+\Big|\int_{-\infty}^{\infty}\int_{0}^{1}\int_{0}^{1}\frac{tf''(x+z+uty)-tf''(x+uty)}{2|y|^{\alpha-1}}
I_{\beta}(y)\dif u\dif t\dif y\Big|\\
&\leq |z|\Big[\int_{-\infty}^{\infty}\frac{\|f''\|_{\infty}}{|y|^{\alpha}}I^{\beta}(y)\dif y+\int_{-\infty}^{\infty}\int_{0}^{1}\frac{t\|f^{(3)}\|_{\infty}}{2|y|^{\alpha-1}}I_{\beta}(y)\dif t\dif y\Big]=\Big(\frac{2\|f''\|_{\infty}}{\alpha-1}+\frac{\|f^{(3)}\|_{\infty}}{2(2-\alpha)}\Big)|z|.
\end{align*}
2. When $\alpha=1,$ we have by (\ref{op})
\begin{align}\label{1}
\frac{1}{d_{\alpha}}\mathcal{L}^{1,\beta}f(x)&=\int_{-\infty}^{\infty}\frac{f(y+x)-f(x)}{2|y|^{2}}I^{\beta}(y)\dif y+\int_{-\infty}^{\infty}\frac{f(y+x)-f(x)-yf'(x)}{2|y|^{2}}I_{\beta}(y)\dif y\nonumber\\
&=\int_{-\infty}^{\infty}\frac{f(y+x)-f(x)}{2|y|^{2}}I^{\beta}(y)\dif y+\int_{-\infty}^{\infty}\int_{0}^{1}\int_{0}^{1}\frac{tf''(x+uty)}{2}
I_{\beta}(y)\dif u\dif t\dif y,
\end{align}
then by the same argument as above, we have
\begin{align*}
&\quad\frac{1}{d_{\alpha}}\big|\mathcal{L}^{1,\beta}f(x+z)-\mathcal{L}^{1,\beta}f(x)\big|\leq\Big(2\|f'\|_{\infty}+\frac{1}{2}\|f^{(3)}\|_{\infty}\Big)|z|.
\end{align*}
3. When $\alpha\in(0,1),$ we have by (\ref{op})
\begin{align}\label{0}
\frac{1}{d_{\alpha}}\mathcal{L}^{\alpha,\beta}f(x)&=\int_{-\infty}^{\infty}\frac{f(y+x)-f(x)}{2|y|^{1+\alpha}}I^{\beta}(y)\dif y+\int_{-\infty}^{\infty}\int_{0}^{1}\frac{yf'(x+ty)}{2|y|^{1+\alpha}}
I_{\beta}(y)\dif t\dif y,
\end{align}
it follows that
\begin{align*}
\frac{1}{d_{\alpha}}\big|\mathcal{L}^{\alpha,\beta}f(x+z)-\mathcal{L}^{\alpha,\beta}f(x)\big|\leq\Big(\frac{2}{\alpha}\|f'\|_{\infty}+\frac{1}{1-\alpha}\|f''\|_{\infty}\Big)|z|.
\end{align*}
\end{proof}

By (\ref{1}) and (\ref{0}) above, we can immediately obtain the following proposition:
\bp\label{p1}
Let $\alpha\in(0,1].$ For any $f\in\mathcal{C}_{b}^{2},$ we have
\begin{align*}
\|\mathcal{L}^{\alpha,\beta}f\|_{\infty}\leq\hat{D}_{\alpha},
\end{align*}
where
\begin{align*}
\hat{D}_{\alpha}=
\begin{cases}
2d_{\alpha}\|f\|_{\infty}+\frac{d_{\alpha}}{2}\|f''\|_{\infty},\quad &\alpha=1,\\
\frac{2d_{\alpha}}{\alpha}\|f\|_{\infty}+\frac{d_{\alpha}}{1-\alpha}\|f'\|_{\infty},\quad &\alpha\in(0,1).
\end{cases}
\end{align*}
\ep

\medskip
In addition, it is easy to verify by the definition of $\mathcal{L}^{\alpha,\beta}$ that if $z=x-y,$ then
\begin{equation}\label{plus}
\mathcal{L}^{\alpha,\beta}_{x}f(x-y)=\mathcal{L}^{\alpha,\beta}_{z}f(z),
\end{equation}
where $\mathcal{L}^{\alpha,\beta}_{x}$ means that the operator $\mathcal{L}^{\alpha,\beta}$ acts on the variable $x$. Similarly, for $z=cx$ for some constant $c>0,$ we have
\begin{equation}\label{time1}
\mathcal{L}^{\alpha,\beta}_{x}f(cx)=c^{\alpha}\mathcal{L}^{\alpha,\beta}_{z}f(z).
\end{equation}
\subsection{Truncation for asymmetric $\alpha-$stable process $\hat{Y}$}
When $\alpha\in(0,1],$ we have by (\ref{G1}) that $\mathbb{E}|\hat{Y}_{s}|=\infty$ for any $s>0$, we need the following lemma for the analysis in the next section.
\bl\label{12}
Consider $\alpha\in(0,1].$ Let $\hat{Y}$ be the one-dimensional asymmetric $\alpha-$stable process, then for any $0<a<1,$ $z\in\mathbb{R}$ and $f\in\mathcal{C}^{3}_{b}(\mathbb{R}),$ we have
\begin{align}\label{11}
\mathbb{E}\big[\int_{0}^{1}\big|\mathcal{L}^{\alpha,\beta}f(z)-\mathcal{L}^{\alpha,\beta}f(a\hat{Y}_{s}+z)\big|\dif s\big]\leq C_{\alpha}
\begin{cases}
(\hat{D}_{\alpha}+D_{\alpha})a-D_{\alpha}a\log a,\quad &\alpha=1,\\
\qquad(\hat{D}_{\alpha}+D_{\alpha})a^{\alpha},\quad &\alpha\in(0,1),
\end{cases}
\end{align}
where $\hat{D}_{\alpha}$ and $D_{\alpha}$ are defined as above.
\el
\begin{proof}
Observe
\begin{align*}
&\quad\mathbb{E}\big[\int_{0}^{1}\big|\mathcal{L}^{\alpha,\beta}f(z)-\mathcal{L}^{\alpha,\beta}f(a\hat{Y}_{s}+z)\big|\dif s\big]\\
&=\mathbb{E}\Big[\int_{0}^{1}\big|\mathcal{L}^{\alpha,\beta}f(z)-\mathcal{L}^{\alpha,\beta}f(a\hat{Y}_{s}+z)\big|\big[{\bf 1}_{(a^{-1},\infty)}(|\hat{Y}_{s}|)+{\bf 1}_{(0,a^{-1})}(|\hat{Y}_{s}|)\big]\dif s\Big]:=\mathcal{I}+\mathcal{II}.
\end{align*}
By Proposition \ref{p1} and Lemma \ref{g1}, we have
\begin{align*}
\mathcal{I}\leq 2\hat{D}_{\alpha}\int_{0}^{1}\mathbb{P}(|\hat{Y}_{s}|\geq a^{-1})\dif s\leq C_{\alpha}\hat{D}_{\alpha}\int_{0}^{1}\int_{a^{-1}}^{\infty}\frac{s}{y^{\alpha+1}}\dif y\dif s\leq C_{\alpha}\hat{D}_{\alpha}a^{\alpha}.
\end{align*}
By Proposition \ref{p} and Lemma \ref{g1}, when $\alpha=1$ and $\beta=0,$ we have
\begin{align*}
\mathcal{II}&\leq D_{\alpha}a\mathbb{E}\big[\int_{0}^{1}|\hat{Y}_{s}|I_{\{|\hat{Y}_{s}|\leq a^{-1}\}}\dif s\big]\\
&\leq C_{\alpha}D_{\alpha}a\int_{0}^{1}\int_{0}^{a^{-1}}ys^{-1}\big(1\wedge\frac{s^{2}}{y^{2}}\big)\dif y\dif s\\
&=C_{\alpha}D_{\alpha}a\int_{0}^{1}\int_{0}^{s}ys^{-1}\dif y\dif s+C_{\alpha}D_{\alpha}a\int_{0}^{1}\int_{s}^{a^{-1}}\frac{s}{y}\dif y\dif s\\
&=C_{\alpha}D_{\alpha}a\Big(\frac{1}{4}+\int_{0}^{1}s(\log a^{-1}+\log s^{-1})\dif s\Big)\\
&\leq C_{\alpha}D_{\alpha}a-C_{\alpha}D_{\alpha}a\log a,
\end{align*}
where the last inequality is by the fact $\int_{0}^{1}s\log s^{-1}\dif s\leq1.$ When $\alpha\in(0,1),$ we have
\begin{align*}
\mathcal{II}&\leq D_{\alpha}a\mathbb{E}\big[\int_{0}^{1}|\hat{Y}_{s}|I_{\{|\hat{Y}_{s}|\leq a^{-1}\}}\dif s\big]\\
&\leq C_{\alpha}D_{\alpha}a\int_{0}^{1}\int_{0}^{a^{-1}}ys^{-\frac{1}{\alpha}}\big(1\wedge\frac{s^{\frac{\alpha+1}{\alpha}}}{y^{\alpha+1}}\big)\dif y\dif s\leq C_{\alpha}D_{\alpha}a^{\alpha}.
\end{align*}
Collecting the previous estimates, we complete the proof.
\end{proof}
\br
In the above lemma, because of $\mathbb{E}|\hat{Y}_{s}|<\infty$ in the case $\alpha\in(1,2),$ we have by Proposition \ref{p} that
\begin{align}\label{22}
\mathbb{E}\big[\int_{0}^{1}\big|\mathcal{L}^{\alpha,\beta}f(z)-\mathcal{L}^{\alpha,\beta}f(a\hat{Y}_{s}+z)\big|\dif s\big]\leq D_{\alpha}a\int_{0}^{1}\mathbb{E}|\hat{Y}_{s}|\dif s\leq C_{\alpha}D_{\alpha}a.
\end{align}
\er

\subsection{Truncation for random variable $X$}
Let $X$ have a distribution of the form (\ref{px}) with $\epsilon(x)$ satisfying (\ref{bound}), then it is obvious that $\mathbb{E}|X|=\infty$ in the case $\alpha\in(0,1].$ However, we can use a truncation technique to handle the problem. Before giving the truncation Lemma, we need
\bl
Let $X\geq0$ be a random variable, for any $t>0$,
\begin{align}\label{expectation}
\mathbb{E}\big[X{\bf 1}_{(0,t)}\big]=\int_{0}^{t}\mathbb{P}(X>r)\dif r-t\mathbb{P}(X>t).
\end{align}
\el
\begin{proof}
Using the definition of expected value, Fubini's theorem and then calculating the resulting integrals gives
\begin{align*}
\int_{0}^{t}\mathbb{P}(X>r)\dif r&=\int_{0}^{t}\int_{\Omega}{\bf 1}_{(r,\infty)}(X)\dif \mathbb{P}\dif r\\
&=\int_{\Omega}\int_{0}^{t}{\bf 1}_{(0,X)}(r)\dif r\dif \mathbb{P}=\int_{\Omega}(X\wedge t)\dif \mathbb{P}=\mathbb{E}\big[X{\bf 1}_{(0,t)}\big]+t\mathbb{P}(X>t),
\end{align*}
from which we immediately obtain the equality in the lemma, as desired.
\end{proof}

Now, we are at the position to give the truncation lemma.
\bl\label{random truncation}
Consider $\alpha\in(0,1].$ Let $X$ have a distribution of the form (\ref{px}) with $\epsilon(x)$ satisfying (\ref{bound}), then for any $0<a<1,$ $z\in\mathbb{R}$ and $f\in\mathcal{C}^{2}_{b}(\mathbb{R}),$ we have
\begin{align*}
\mathbb{E}|f'(z+aX)-f'(z)|\leq
\begin{cases}
2\|f''\|_{\infty}a+2(A+K)\big(2\|f'\|_{\infty}+\|f''\|_{\infty}\log a^{-1}\big)a, &\alpha=1,\\
2\|f''\|_{\infty}a+2(A+K)\big(2\|f'\|_{\infty}+\frac{\|f''\|_{\infty}}{1-\alpha}\big)a^{\alpha}, &\alpha\in(0,1).
\end{cases}
\end{align*}
\el
\begin{proof}
Observe
\begin{align*}
\mathbb{E}|f'(z+aX)-f'(z)|&=\mathbb{E}\Big[|f'(z+aX)-f'(z)|\big({\bf 1}_{(a^{-1},\infty)}(|X|)+{\bf 1}_{(0,a^{-1}]}(|X|)\big)\Big]\\
&\leq2\|f'\|_{\infty}\mathbb{P}\big(|X|>a^{-1}\big)+\|f''\|_{\infty}a\mathbb{E}\big[|X|{\bf 1}_{(0,a^{-1}]}(|X|)\big].
\end{align*}
By (\ref{bound}), we know $|\epsilon(x)|\leq K$ for $|x|\geq1,$ so
\begin{align*}
\mathbb{P}\big(|X|>a^{-1}\big)=a^{\alpha}\big(A+\epsilon(a^{-1})\big)(1+\beta)+a^{\alpha}\big(A+\epsilon(-a^{-1})\big)(1-\beta)\leq2(A+K)a^{\alpha}.
\end{align*}
By (\ref{expectation}),
\begin{align*}
\mathbb{E}\big[|X|{\bf 1}_{(0,a^{-1}]}(|X|)\big]&\leq\int_{0}^{a^{-1}}\mathbb{P}(|X|>r)\dif r\\
& \le 1+(1+\beta)\int_{1}^{a^{-1}}\frac{A+\epsilon(x)}{x^{\alpha}}\dif x+(1-\beta)\int_{1}^{a^{-1}}\frac{A+\epsilon(-x)}{x^{\alpha}}\dif x,
\end{align*}
1.) when $\alpha=1,$ we have
\begin{align*}
(1+\beta)\int_{1}^{a^{-1}}\frac{A+\epsilon(x)}{x}\dif x\leq (1+\beta)\int_{1}^{a^{-1}}\frac{A+K}{x}\dif x=(1+\beta)(A+K)\log a^{-1}.
\end{align*}
2.) when $\alpha\in(0,1),$ we have
\begin{align*}
(1+\beta)\int_{1}^{a^{-1}}\frac{A+\epsilon(x)}{x^{\alpha}}\dif x\leq (1+\beta)\int_{1}^{a^{-1}}\frac{A+K}{x^{\alpha}}\dif x\leq \frac{(1+\beta)(A+K)}{1-\alpha}a^{\alpha-1}.
\end{align*}
Similar bounds hold true for $(1-\beta)\int_{0}^{a^{-1}}\frac{A+\epsilon(-x)}{x^{\alpha}}\dif x$. Collecting the above estimates, we conclude the proof.
\end{proof}

\br
From the proof of Lemma \ref{random truncation}, we immediately have
\begin{align}\label{te}
\mathbb{E}\big[|X|{\bf 1}_{(0,a^{-1}]}(|X|)\big]\leq2
\begin{cases}
1+(A+K)\log a^{-1}, &\alpha=1,\\
1+\frac{A+K}{1-\alpha}a^{\alpha-1}, &\alpha\in(0,1).
\end{cases}
\end{align}
\er
\subsection{Taylor-like expansions for Theorem \ref{main}}
In this section, we develop the following Taylor-like expansions, which can be taken as replacements of the Taylor expansions in the Lindeberg's approach to proving the normal CLT.

\vskip 3mm

$\bullet$ \underline{{\bf $\alpha\in(1,2)$}} {\bf :}
\bl\label{ml}
Consider $\alpha\in(1,2).$ Let $X$ have a distribution $F_{X}$ with the form (\ref{px}), $X$ and $Z$ are independent. For any $0<a\leq(2A)^{-\frac{1}{\alpha}}\wedge1$ and $f\in\mathcal{C}^{2}_{b}(\mathbb{R})$, we have
\begin{align*}
&\quad\Big|\mathbb{E}\big[f(Z+aX)\big]-\mathbb{E}[f(Z)]-\mathbb{E}[aX]\mathbb{E}[f'(Z)]-\frac{2A\alpha}{d_{\alpha}}a^{\alpha}\mathbb{E}[\mathcal{L}^{\alpha,\beta}f(Z)]\Big|\\
&\leq\frac{4\|f''\|_{\infty}}{2-\alpha}(2A)^{\frac{2}{\alpha}}a^{2}+\frac{8\|f'\|_{\infty}}{\alpha-1}a^{\alpha}\sup_{|x|\geq a^{-1}}|\epsilon(x)|+2\|f''\|_{\infty}a^{2}\int_{-a^{-1}}^{a^{-1}}\frac{|\epsilon(x)|}{x^{\alpha-1}}\dif x.
\end{align*}
\el
\begin{proof}
Denote ${\bf 1}_{\beta}(y)=(1+\beta){\bf 1}_{(0,\infty)}(y)+(1-\beta){\bf 1}_{(-\infty,0)}(y)$, we have by (\ref{op})
\begin{align*}
\frac{2A\alpha}{d_{\alpha}}a^{\alpha}\mathbb{E}[\mathcal{L}^{\alpha,\beta}f(Z)]&=A\alpha a^{\alpha}\,
\mathbb{E}\Big[\int_\R \frac{f(Z+y)-f(Z)-yf'(Z)}{|y|^{1+\alpha}}{\bf 1}_{\beta}(y)
 dy\Big]\\
 &=A\alpha\,\mathbb{E}
\Big[\int_\R \frac{f(Z+ax)-f(Z)-axf'(Z)}{|x|^{1+\alpha}}
{\bf 1}_{\beta}(x)
 dx\Big]\\
 &=A\alpha\,\mathbb{E}
\Big[\int_{|x|\geq(2A)^{\frac{1}{\alpha}}} \frac{f(Z+ax)-f(Z)-axf'(Z)}{|x|^{1+\alpha}}
{\bf 1}_{\beta}(x)
 dx\Big]+\mathcal{R},
\end{align*}
where the second equality is by taking $y=ax$ and
\begin{align}\label{R2}
\mathcal{R}=A\alpha\,\mathbb{E}
\Big[\int_{|x|<(2A)^{\frac{1}{\alpha}}} \frac{f(Z+ax)-f(Z)-axf'(Z)}{|x|^{1+\alpha}}
{\bf 1}_{\beta}(x)
 dx\Big].
\end{align}
Since $A\alpha\int_{|x|\geq(2A)^{\frac{1}{\alpha}}} \frac{1}{|x|^{1+\alpha}}
{\bf 1}_{\beta}(x)
 dx=1,$ we can consider a random variable $\tilde{X}$ which is independent of $Z$ and satisfies
\begin{align}\label{F}
\mathbb{P}(\tilde{X}>x)=\frac{A(1+\beta)}{|x|^{\alpha}},\quad x\geq(2A)^{\frac{1}{\alpha}},\qquad \mathbb{P}(\tilde{X}\leq x)=\frac{A(1-\beta)}{|x|^{\alpha}},\quad x\leq-(2A)^{\frac{1}{\alpha}},
\end{align}
it follows that
 \begin{align*}
 \frac{2A\alpha}{d_{\alpha}}a^{\alpha}\mathbb{E}[\mathcal{L}^{\alpha,\beta}f(Z)]=\mathbb{E}\big[f(Z+a\tilde{X})-f(Z)-a\tilde{X}f'(Z)\big]+\mathcal{R}.
 \end{align*}
As a result, denote the distribution function of $\tilde{X}$ by $F_{\tilde{X}},$ then
\begin{align}\label{1result}
&\quad\Big|\mathbb{E}\big[f(Z+aX)\big]-\mathbb{E}[f(Z)]-\mathbb{E}[aX]\mathbb{E}[f'(Z)]-\frac{2A\alpha}{d_{\alpha}}a^{\alpha}\mathbb{E}[\mathcal{L}^{\alpha,\beta}f(Z)]\Big|\nonumber\\
&\leq\mathbb{E}\Big|\int_{-\infty}^{\infty}\big[f(Z+ax)-axf'(Z)\big]\dif \big[F_{X}(x)-F_{\tilde{X}}(x)\big]\Big|+|\mathcal{R}|.
\end{align}
By (\ref{px}) and (\ref{F}), we have
\begin{align*}
F_{X}(x)-F_{\tilde{X}}(x)=&\big(\frac{1}{2}-\frac{A+\epsilon(x)}{|x|^{\alpha}}\big)(1+\beta){\bf 1}_{(0,(2A)^{\frac{1}{\alpha}})}(x)-\frac{\epsilon(x)}{|x|^{\alpha}}(1+\beta){\bf 1}_{((2A)^{\frac{1}{\alpha}},\infty)}(x)\\
&+\big(\frac{A+\epsilon(x)}{|x|^{\alpha}}-\frac{1}{2}\big)(1-\beta){\bf 1}_{(-(2A)^{\frac{1}{\alpha}},0)}(x)+\frac{\epsilon(x)}{|x|^{\alpha}}(1-\beta){\bf 1}_{(-\infty,-(2A)^{\frac{1}{\alpha}})}(x),
\end{align*}
using integration by parts, we have
\begin{equation} \label{e:IBP}
\begin{split}
&\mathbb{E}\Big|\int_{-\infty}^{\infty}\big[f(Z+ax)-axf'(Z)\big]\dif \big[F_{X}(x)-F_{\tilde{X}}(x)\big]\Big|\\
=&\mathbb{E}\Big|\int_{-\infty}^{\infty}\big[F_{X}(x)-F_{\tilde{X}}(x)\big]\big[af'(Z+ax)-af'(Z)\big]\dif x\Big|\\
\leq&2\mathbb{E}\Big[\int_{-(2A)^{\frac{1}{\alpha}}}^{(2A)^{\frac{1}{\alpha}}}\big|af'(Z+ax)-af'(Z)\big|\dif x\Big]+2\mathbb{E}\Big|\int_{(2A)^{\frac{1}{\alpha}}}^{\infty}\big[af'(Z+ax)-af'(Z)\big]\frac{\epsilon(x)}{x^{\alpha}}\dif x\Big|\\
&+2\mathbb{E}\Big|\int_{(2A)^{\frac{1}{\alpha}}}^{\infty}\big[af'(Z-ax)-af'(Z)\big]\frac{\epsilon(-x)}{x^{\alpha}}\dif x\Big|,
\end{split}
\end{equation}
and
\begin{align*}
2\mathbb{E}\Big[\int_{-(2A)^{\frac{1}{\alpha}}}^{(2A)^{\frac{1}{\alpha}}}\big|af'(Z+ax)-af'(Z)\big|\dif x\Big]\leq2(2A)^{\frac{2}{\alpha}}\|f''\|_{\infty}a^{2}.
\end{align*}
For the remainder, one has
\begin{align*}
\mathbb{E}\Big[\int_{a^{-1}}^{\infty}\big|af'(Z+ax)-af'(Z)\big|\frac{|\epsilon(x)|}{x^{\alpha}}\dif x\Big]\leq\frac{2\|f'\|_{\infty}}{\alpha-1}a^{\alpha}\sup_{x\geq a^{-1}}|\epsilon(x)|,
\end{align*}
whereas
\begin{align*}
\mathbb{E}\Big[\int_{(2A)^{\frac{1}{\alpha}}}^{a^{-1}}\big|af'(Z+ax)-af'(Z)\big|\frac{|\epsilon(x)|}{x^{\alpha}}\dif x\Big]\leq\|f''\|_{\infty}a^{2}\int_{0}^{a^{-1}}\frac{|\epsilon(x)|}{x^{\alpha-1}}\dif x.
\end{align*}
Since similar bounds hold true for $\mathbb{E}\Big[\int_{(2A)^{\frac{1}{\alpha}}}^{\infty}\big|af'(Z-ax)-af'(Z)\big|\frac{\epsilon(-x)}{x^{\alpha}}\dif x\Big]$ and
\begin{align}\label{R1}
|\mathcal{R}|\leq 2A\alpha\|f''\|_{\infty}a^{2}\int_{|x|<(2A)^{\frac{1}{\alpha}}} \frac{1}{|x|^{\alpha-1}}dx=\frac{4A\alpha\|f''\|_{\infty}}{2-\alpha}(2A)^{\frac{2-\alpha}{\alpha}}a^{2}.
\end{align}
the desired conclusion follows.
\end{proof}

$\bullet$ \underline{{\bf $\alpha=1$ and $\beta=0$}} {\bf :}

\bl\label{ml1}
Consider $\alpha=1,$ $\beta=0$ and $\gamma\in(0,\infty).$ Let $X$ have a distribution of the form (\ref{px}) with $\epsilon(x)$ satisfying (\ref{bound}), $X$ and $Z$ are independent. We have, for any $0<a\leq(2A)^{-1}\wedge1$ and $f\in\mathcal{C}^{2}_{b}(\mathbb{R}),$ denote
\begin{align*}
T_{1}:=\Big|\mathbb{E}\big[f(Z+aX)\big]-\mathbb{E}[f(Z)]-\mathbb{E}[aX{\bf 1}_{(-1,1)}(aX)]\mathbb{E}[f'(Z)]-\frac{2A}{d_{1}}a\mathbb{E}[\mathcal{L}^{1,0}f(Z)]\Big|,
\end{align*}
then we have
\begin{align*}
T_{1}\leq12A^{2}\|f''\|_{\infty}a^{2}&+2(2\|f\|_{\infty}+\|f'\|_{\infty})a\sup_{|x|\geq a^{-1}}|\epsilon(x)|\\
&+\|f''\|_{\infty}a^{2}\int_{-a^{-1}}^{a^{-1}}|\epsilon(x)|dx+\|f'\|_{\infty}a\int_{a^{-1}}^{\infty}\frac{|\epsilon(x)|+|\epsilon(-x)|}{x}dx.
\end{align*}
\el
\begin{proof}
By the same argument as (\ref{1result}), we have
\begin{align}\label{2result}
T_{1}\leq\mathbb{E}\Big|\int_{-\infty}^{\infty}\big[f(Z+ax)-ax{\bf 1}_{(-1,1)}(ax)f'(Z)\big]\dif \big[F_{X}(x)-F_{\tilde{X}}(x)\big]\Big|+|\mathcal{R}|,
\end{align}
where $F_{\tilde{X}}$ and $\mathcal{R}$ is defined by (\ref{F}) and (\ref{R2}) with $\alpha=1,\beta=0,$ respectively. What's more, by (\ref{R1}), we know $|\mathcal{R}|\leq8A^{2}\|f''\|_{\infty}a^{2}.$\\
For the first term, using an integration by parts similar to \eqref{e:IBP} and (\ref{bound}), we have
\begin{align*}
\mathbb{E}\Big|\int_{a^{-1}}^{\infty}f(Z+ax)\dif \big[F_{X}(x)-F_{\tilde{X}}(x)\big]\Big|&\leq\|f\|_{\infty}a|\epsilon(a^{-1})|+\|f'\|_{\infty}a\int_{a^{-1}}^{\infty}\frac{|\epsilon(x)|}{x}\dif x
\end{align*}
and in the same way
\begin{align*}
&\mathbb{E}\Big|\int^{-a^{-1}}_{-\infty}f(Z+ax)\dif \big[F_{X}(x)-F_{\tilde{X}}(x)\big]\Big|\leq \|f\|_{\infty}a|\epsilon(-a^{-1})|+\|f'\|_{\infty}a\int^{-a^{-1}}_{-\infty}\frac{|\epsilon(x)|}{|x|}\dif x,
\end{align*}
whereas
\begin{align*}
&\mathbb{E}\Big|\int_{-a^{-1}}^{a^{-1}}\big[f(Z+ax)-axf'(Z)\big]\dif \big[F_{X}(x)-F_{\tilde{X}}(x)\big]\Big|\\
\leq&\big(\|f\|_{\infty}+\|f'\|_{\infty}\big)a\big(|\epsilon(a^{-1})|+|\epsilon(-a^{-1})|\big)+\|f''\|_{\infty}a^{2}\big(\int_{-2A}^{2A}|x|\dif x+\int_{-a^{-1}}^{a^{-1}}|\epsilon(x)|{\bf 1}_{(2A,\infty)}(|x|)\dif x\big)\\
\leq&2\big(\|f\|_{\infty}+\|f'\|_{\infty}\big)a\sup_{|x|\geq a^{-1}}|\epsilon(x)|+4A^{2}\|f''\|_{\infty}a^{2}+\|f''\|_{\infty}a^{2}\int_{-a^{-1}}^{a^{-1}}|\epsilon(x)|\dif x,
\end{align*}
 the desired conclusion follows.
\end{proof}

$\bullet$ \underline{{\bf $\alpha\in(0,1)$}} {\bf :} For any $\beta\in[-1,1],$ we have
\begin{align*}
\int_{\mathbb{R}}\frac{y{\bf 1}_{(-1,1)}(y)}{2|y|^{1+\alpha}}\left[
(1+\beta){\bf 1}_{(0,\infty)}(y)+(1-\beta){\bf 1}_{(-\infty,0)}(y)
\right]
 dy=\frac{\beta}{1-\alpha},
\end{align*}
which follows that
\begin{align}\label{ashift}
&\frac{1}{d_{\alpha}}\mathcal{L}^{\alpha,\beta}f(x)-\frac{\beta f'(x)}{1-\alpha}\nonumber\\
=&\int_\R \frac{f(y+x)-f(x)-y{\bf 1}_{(-1,1)}(y)f'(x)}{2|y|^{1+\alpha}}
\left[
(1+\beta){\bf 1}_{(0,\infty)}(y)+(1-\beta){\bf 1}_{(-\infty,0)}(y)
\right]
 dy.
\end{align}
According to (\ref{ashift}), we have the following Taylor-like expansion lemma.
\bl\label{ml2}
Consider $\alpha\in(0,1).$ Let $X$ have a distribution $F_{X}$ with the form (\ref{px}) satisfying (\ref{bound}), $X$ and $Z$ are independent. We have, for any $0<a\leq(2A)^{-\frac{1}{\alpha}}\wedge1$ and $f\in\mathcal{C}^{2}_{b}(\mathbb{R}),$ denote
\begin{align*}
T_{2}\!:=\!\Big|\mathbb{E}\big[f(Z+aX)\big]\!-\!\mathbb{E}[f(Z)]\!-\!\mathbb{E}[aX{\bf 1}_{(-1,1)}(aX)]\mathbb{E}[f'(Z)]-\frac{2A\alpha}{d_{\alpha}}a^{\alpha}\mathbb{E}[\mathcal{L}^{\alpha,\beta}f(Z)-\frac{\beta d_{\alpha}f'(Z)}{1-\alpha}]\Big|.
\end{align*}
a.) When $\gamma\in(1-\alpha,\infty),$ we have
\begin{align*}
T_{2}\leq\frac{2+\alpha}{2-\alpha}(2A)^{\frac{2}{\alpha}}\|f''\|_{\infty}a^{2}&+2(3\|f\|_{\infty}
+2\|f'\|_{\infty})a^{\alpha}\sup_{|x|\geq a^{-1}}|\epsilon(x)|\\&+2\|f''\|_{\infty}a^{2}\int_{-a^{-1}}^{a^{-1}}\frac{|\epsilon(x)|}{|x|^{\alpha-1}}dx+\|f'\|_{\infty}a\int_{a^{-1}}^{\infty}\frac{|\epsilon(x)|+|\epsilon(-x)|}{x^{\alpha}}dx.
\end{align*}
b.) When $\gamma\in[0,1-\alpha],$ we have
\begin{align*}
T_{2}\leq\frac{2+\alpha}{2-\alpha}(2A)^{\frac{2}{\alpha}}\|f''\|_{\infty}a^{2}&+\big[(4A+6K+4)\|f\|_{\infty}+\frac{8-4\alpha}{1-\alpha}\|f'\|_{\infty}\big]a^{\alpha}\sup_{|x|\geq a^{-1}}|\epsilon(x)|^{\alpha}\\
&+2\|f''\|_{\infty}a^{2}\int_{-a^{-1}}^{a^{-1}}\frac{|\epsilon(x)|}{|x|^{\alpha-1}}dx.
\end{align*}
\el
\begin{proof}
By the same argument as (\ref{1result}), we have
\begin{align}\label{2result}
T_{2}\leq\mathbb{E}\Big|\int_{-\infty}^{\infty}\big[f(Z+ax)-ax{\bf 1}_{(-1,1)}(ax)f'(Z)\big]\dif \big[F_{X}(x)-F_{\tilde{X}}(x)\big]\Big|+|\mathcal{R}|,
\end{align}
where $F_{\tilde{X}}$ and $\mathcal{R}$ is defined by (\ref{F}) and (\ref{R2}) with $\alpha\in(0,1),$ respectively. What's more, by (\ref{R1}), we know $|\mathcal{R}|\leq\frac{4A\alpha\|f''\|_{\infty}}{2-\alpha}(2A)^{\frac{2-\alpha}{\alpha}}a^{2}.$\\
For the first term, according to (\ref{bound}),\\
1. When $\gamma\in(1-\alpha,\infty),$ using an integration by parts similar to \eqref{e:IBP} and (\ref{bound}), we have
\begin{align*}
\mathbb{E}\Big|\int_{a^{-1}}^{\infty}f(Z+ax)\dif \big[F_{X}(x)-F_{\tilde{X}}(x)\big]\Big|&\leq\|f\|_{\infty}a^{\alpha}|\epsilon(a^{-1})|+\|f'\|_{\infty}a\int_{a^{-1}}^{\infty}\frac{|\epsilon(x)|}{x^{\alpha}}\dif x.
\end{align*}
Similarly,
\begin{align*}
&\mathbb{E}\Big|\int^{-a^{-1}}_{-\infty}f(Z+ax)\dif \big[F_{X}(x)-F_{\tilde{X}}(x)\big]\Big|\leq \|f\|_{\infty}a^{\alpha}|\epsilon(-a^{-1})|+\|f'\|_{\infty}a\int^{-a^{-1}}_{-\infty}\frac{|\epsilon(x)|}{|x|^{\alpha}}\dif x.
\end{align*}
Moreover, using integration by parts again, we obtain
\begin{align*}
&\mathbb{E}\Big|\int_{-a^{-1}}^{a^{-1}}\big[f(Z+ax)-axf'(Z)\big]\dif \big[F_{X}(x)-F_{\tilde{X}}(x)\big]\Big|\\
\leq&4\big(\|f\|_{\infty}+\|f'\|_{\infty}\big)a^{\alpha}\sup_{|x|\geq a^{-1}}|\epsilon(x)|+\|f''\|_{\infty}a^{2}\Big(\int_{-(2A)^{\frac{1}{\alpha}}}^{(2A)^{\frac{1}{\alpha}}}|x|\dif x+2\int_{-a^{-1}}^{a^{-1}}\frac{|\epsilon(x)|}{|x|^{\alpha-1}}\dif x\Big)\\
\leq&4\big(\|f\|_{\infty}+\|f'\|_{\infty}\big)a^{\alpha}\sup_{|x|\geq a^{-1}}|\epsilon(x)|+(2A)^{\frac{2}{\alpha}}\|f''\|_{\infty}a^{2}+2\|f''\|_{\infty}a^{2}\int_{-a^{-1}}^{a^{-1}}\frac{|\epsilon(x)|}{|x|^{\alpha-1}}\dif x.
\end{align*}
2. When $\gamma\in[0,1-\alpha],$ we choose a number $N>a^{-1}.$ One has by $|\epsilon(x)|\leq K$ for $|x|>N,$
\begin{align*}
\mathbb{E}\Big|\int_{|x|>N}f(Z+ax)\dif \big[F_{X}(x)-F_{\tilde{X}}(x)\big]\Big|&\leq\|f\|_{\infty}\big[\int_{|x|>N}\dif F_{X}(x)+\int_{|x|>N}\dif F_{\tilde{X}}(x)\big]\\
&\leq(4A+2K)\|f\|_{\infty}N^{-\alpha},
\end{align*}
whereas by integration by parts
\begin{align*}
&\mathbb{E}\Big|\big(\int_{-a^{-1}}^{a^{-1}}+\int_{a^{-1}<|x|<N}\big)\big[f(Z+ax)-ax{\bf 1}_{(-1,1)}(ax)f'(Z)\big]\dif \big[F_{X}(x)-F_{\tilde{X}}(x)\big]\Big|\\
%\leq&2(\|f\|_{\infty}\!+\!\|f'\|_{\infty})a^{\alpha}\big(|\epsilon(a^{-1})|\!+\!|\epsilon(-a^{-1})|\big)\!\!+\!\mathbb{E}\Big[\int_{-a^{-1}}^{a^{-1}}\!\!\big|F_{X}(x)-F_{\tilde{X}}(x)\big||af'(Z+ax)-af'(Z)|\dif x\Big]\\
%&+2\|f\|_{\infty}N^{-\alpha}\big(|\epsilon(N)|+|\epsilon(-N)|\big)+\mathbb{E}\Big[\int_{a^{-1}\leq|X|\leq N}\big|F_{X}(x)-F_{\tilde{X}}(x)\big|a|f'(Z+ax)|\dif x\Big]\\
\leq&a^{2}\|f''\|_{\infty}\big((2A)^{\frac{2}{\alpha}}+2\int_{-a^{-1}}^{a^{-1}}\frac{|\epsilon(x)|}{|x|^{\alpha-1}}\dif x\big)+4(\|f\|_{\infty}+\|f'\|_{\infty})a^{\alpha}\sup_{|x|\geq a^{-1}}|\epsilon(x)|\\
&+4K\|f\|_{\infty}N^{-\alpha}+2\|f'\|_{\infty}a\Big(\int_{a^{-1}}^{N}\frac{|\epsilon(x)|}{x^{\alpha}}\dif x+\int_{a^{-1}}^{N}\frac{|\epsilon(-x)|}{x^{\alpha}}\dif x\Big).
\end{align*}
Since
\begin{align*}
\int_{a^{-1}}^{N}\frac{|\epsilon(x)|}{x^{\alpha}}\dif x+\int_{a^{-1}}^{N}\frac{|\epsilon(-x)|}{x^{\alpha}}\dif x\leq\frac{2}{1-\alpha}\sup_{|x|\geq a^{-1}}|\epsilon(x)|N^{1-\alpha},
\end{align*}
we can consider
\begin{align*}
N^{-\alpha}=aN^{1-\alpha}\sup_{|x|\geq a^{-1}}|\epsilon(x)|,
\end{align*}
which implies
\begin{align*}
N=a^{-1}\big(\sup_{|x|\geq a^{-1}}|\epsilon(x)|\big)^{-1},
\end{align*}
the desired conclusion follows.
\end{proof}
~\\

\section{Proof of Theorem \ref{main} and Theorem \ref{improve}}\label{proof}
In this section, with the help of the Taylor-like expansion in the previous section, we prove the main results by the Lindeberg principle.
\begin{proof} [{\bf Proof of Theorem \ref{main}}]
Recall $Y\sim S_{\alpha}(1,\beta)$ in Theorem \ref{main}, let $Y_{1},Y_{2},\cdots$ be independent copies of $Y,$ it is clear $\frac{Y_{1}+\cdots+Y_{n}}{n^{\frac{1}{\alpha}}}\sim S_{\alpha}(1,\beta).$ Recall the definition of $S_{n}$, we denote
\begin{align*}
\hat{X}_{i}=\frac{1}{\sigma}
\begin{cases}
X_{i}-\mathbb{E}[X_{i}], &\alpha\in(1,2),\\
X_{i}-\mathbb{E}\big[X_{i}{\bf 1}_{(0,\sigma n^{\frac{1}{\alpha}})}(|X_{i}|)\big], &\alpha=1,\\
X_{i},  &\alpha\in(0,1),
\end{cases}
\end{align*}
where $i=1,2,\cdots,$ then for any fixed $n$. Set  $Z_{1}=\hat{X}_{2}+\cdots+\hat{X}_{n}$,
\begin{align*}
Z_{i}=Y_{1}+\cdots+Y_{i-1}+\hat{X}_{i+1}+\cdots+\hat{X}_{n},\quad 1<i<n,
\end{align*}
and $Z_{n}=Y_{1}+\cdots+Y_{n-1}$, we have
\begin{align}\label{aL}
\mathbb{E}[f(S_{n})]-\mathbb{E}[f(Y)]&=\mathbb{E}\Big[f\Big(\frac{\hat{X}_{1}+\cdots+\hat{X}_{n}}{n^{\frac{1}{\alpha}}}\Big)\Big]-\mathbb{E}\Big[f\Big(\frac{Y_{1}+\cdots+Y_{n}}{n^{\frac{1}{\alpha}}}\Big)\Big]\nonumber\\
&=\sum_{i=1}^{n}\mathbb{E}\Big[f\Big(\frac{\hat{X}_{i}+Z_{i}}{n^{\frac{1}{\alpha}}}\Big)-f\Big(\frac{Y_{i}+Z_{i}}{n^{\frac{1}{\alpha}}}\Big)\Big]\nonumber\\
&=\sum_{i=1}^{n}\mathbb{E}\Big[f\Big(\frac{\hat{X}_{i}+Z_{i}}{n^{\frac{1}{\alpha}}}\Big)-f\Big(\frac{Z_{i}}{n^{\frac{1}{\alpha}}}\Big)\Big]-
\sum_{i=1}^{n}\mathbb{E}\Big[f\Big(\frac{Y_{i}+Z_{i}}{n^{\frac{1}{\alpha}}}\Big)-f\Big(\frac{Z_{i}}{n^{\frac{1}{\alpha}}}\Big)\Big].
\end{align}
For the first term, since $X_{i}$ and $Z_{i}$ are independent, we have by \eqref{px}
\begin{align*}
\sum_{i=1}^{n}\mathbb{E}\Big[f\Big(\frac{\hat{X}_{i}+Z_{i}}{n^{\frac{1}{\alpha}}}\Big)-f\Big(\frac{Z_{i}}{n^{\frac{1}{\alpha}}}\Big)\Big]
=\frac{1}{n}\sum_{i=1}^{n} \mathbb E \Big[\mathcal{L}^{\alpha,\beta}f\Big(\frac{Z_{i}}{n^{\frac{1}{\alpha}}}\Big)\Big]+\mathbf{I}+\mathbf{II},
\end{align*}
where in the case $\alpha\in[1,2),$
\begin{align*}
\mathbf{I}=\sum_{i=1}^{n}\Big\{\mathbb{E}\Big[f\Big(\frac{Z_{i}}{n^{\frac{1}{\alpha}}}
+\frac{n^{-\frac{1}{\alpha}}}{\sigma}X_{i}\Big)\Big]-\mathbb{E}\Big[f\Big(\frac{Z_{i}}{n^{\frac{1}{\alpha}}}\Big)\Big]
&-\mathbb{E}\Big[\big(\frac{n^{-\frac{1}{\alpha}}}{\sigma}X_{i}\big)^{(\alpha)}\Big]\mathbb{E}\Big[f'\Big(\frac{Z_{i}}{n^{\frac{1}{\alpha}}}\Big)\Big]\\
&-\frac{2A\alpha}{d_{\alpha}}\Big(\frac{n^{-\frac{1}{\alpha}}}{\sigma}\Big)^{\alpha}\mathbb{E}\Big[\mathcal{L}^{\alpha,\beta}f\Big(\frac{Z_{i}}{n^{\frac{1}{\alpha}}}\Big)\Big]\Big\},
\end{align*}
\begin{align*}
\mathbf{II}=\sum_{i=1}^{n}\Big\{\mathbb{E}\Big[f\Big(\frac{Z_{i}}{n^{\frac{1}{\alpha}}}+\frac{X_{i}}{\sigma n^{\frac{1}{\alpha}}}-\mathbb{E}\big[\big(\frac{n^{-\frac{1}{\alpha}}}{\sigma}X_{i}\big)^{(\alpha)}\big]\Big)\Big]&-\mathbb{E}\Big[f\Big(\frac{Z_{i}}{n^{\frac{1}{\alpha}}}+\frac{X_{i}}{\sigma n^{\frac{1}{\alpha}}}\Big)\Big]\\
&+\mathbb{E}\Big[\big(\frac{n^{-\frac{1}{\alpha}}}{\sigma}X_{i}\big)^{(\alpha)}\Big]\mathbb{E}\Big[f'\Big(\frac{Z_{i}}{n^{\frac{1}{\alpha}}}\Big)\Big]\Big\},
\end{align*}
with
\begin{align*}
\big(\frac{n^{-\frac{1}{\alpha}}}{\sigma}X_{i}\big)^{(\alpha)}=
\begin{cases}
\frac{n^{-\frac{1}{\alpha}}}{\sigma}X_{i}, &\alpha\in(1,2),\\
\frac{n^{-\frac{1}{\alpha}}}{\sigma}X_{i}{\bf 1}_{(0,\sigma n^{\frac{1}{\alpha}})}(|X_{i}|), &\alpha=1,
\end{cases}
\end{align*}
and when $\alpha\in(0,1),$ we have
\begin{align*}
\mathbf{I}=\sum_{i=1}^{n}\Big\{\mathbb{E}\Big[f\Big(\frac{Z_{i}}{n^{\frac{1}{\alpha}}}
+\frac{n^{-\frac{1}{\alpha}}}{\sigma}X_{i}\Big)\Big]&-\mathbb{E}\Big[f\Big(\frac{Z_{i}}{n^{\frac{1}{\alpha}}}\Big)\Big]
-\mathbb{E}\Big[\frac{n^{-\frac{1}{\alpha}}}{\sigma}X_{i}{\bf 1}_{(0,\sigma n^{\frac{1}{\alpha}})}(|X_{i}|)\Big]\mathbb{E}\Big[f'\Big(\frac{Z_{i}}{n^{\frac{1}{\alpha}}}\Big)\Big]\\
&-\frac{2A\alpha}{d_{\alpha}}\Big(\frac{n^{-\frac{1}{\alpha}}}{\sigma}\Big)^{\alpha}\mathbb{E}\Big[\mathcal{L}^{\alpha,\beta}f\Big(\frac{Z_{i}}{n^{\frac{1}{\alpha}}}\Big)-\frac{\beta d_{\alpha}}{1-\alpha}f'\Big(\frac{Z_{i}}{n^{\frac{1}{\alpha}}}\Big)\Big]\Big\},
\end{align*}
\begin{align*}
\mathbf{II}=\sum_{i=1}^{n}\Big\{\mathbb{E}\Big[\frac{n^{-\frac{1}{\alpha}}}{\sigma}X_{i}{\bf 1}_{(0,\sigma n^{\frac{1}{\alpha}})}(|X_{i}|)\Big]\mathbb{E}\Big[f'\Big(\frac{Z_{i}}{n^{\frac{1}{\alpha}}}\Big)\Big]-\frac{2A\alpha \beta}{1-\alpha}\Big(\frac{n^{-\frac{1}{\alpha}}}{\sigma}\Big)^{\alpha}\mathbb{E}\Big[f'\Big(\frac{Z_{i}}{n^{\frac{1}{\alpha}}}\Big)\Big]\Big\}.
\end{align*}
For the second term, notice that $Y_{i}$ has the same distribution as $\hat{Y}_{1}$, where $(\hat{Y}_{t})_{t \ge 0}$ is an $\alpha-$stable process determined by \eqref{e:ChfY}, then we have by (\ref{e:SpDyk})
\begin{align*}
\sum_{i=1}^{n}\mathbb{E}\Big[f\Big(\frac{Y_{i}+Z_{i}}{n^{\frac{1}{\alpha}}}\Big)-f\Big(\frac{Z_{i}}{n^{\frac{1}{\alpha}}}\Big)\Big]
&=\sum_{i=1}^{n}\int_{0}^{1}\mathbb{E}\Big[\mathcal{L}^{\alpha,\beta}_{\hat{Y}_{s}}f\Big(\frac{Z_{i}+\hat{Y}_{s}}{n^{\frac{1}{\alpha}}}\Big)\Big]\dif s\\
&=\frac{1}{n}\sum_{i=1}^{n}\int_{0}^{1}\mathbb{E}\Big[\mathcal{L}^{\alpha,\beta}f\Big(\frac{Z_{i}+\hat{Y}_{s}}{n^{\frac{1}{\alpha}}}\Big)\Big]\dif s,
\end{align*}
where the second equality thanks to (\ref{plus}) and (\ref{time1}).\\
Therefore, we have
\begin{align*}
\mathbb{E}[f(S_{n})]-\mathbb{E}[f(Y)]&=\frac{1}{n}\sum_{i=1}^{n}\int_{0}^{1}\mathbb{E}\Big[\mathcal{L}^{\alpha,\beta}f\Big(\frac{Z_{i}}{n^{\frac{1}{\alpha}}}\Big)-\mathcal{L}^{\alpha,\beta}
f\Big(\frac{Z_{i}+\hat{Y}_{s}}{n^{\frac{1}{\alpha}}}\Big)\Big]\dif s+\mathbf{I}+\mathbf{II},
\end{align*}
and using (\ref{22}) and (\ref{11}) respectively, we have
\begin{align}\label{log1}
\frac{1}{n}\sum_{i=1}^{n}\int_{0}^{1}\mathbb{E}\Big|\mathcal{L}^{\alpha,\beta}f\Big(\frac{Z_{i}}{n^{\frac{1}{\alpha}}}\Big)-\mathcal{L}^{\alpha,\beta}
f\Big(\frac{Z_{i}+\hat{Y}_{s}}{n^{\frac{1}{\alpha}}}\Big)\Big|\dif s\leq C_{\alpha}
\begin{cases}
D_{\alpha}n^{-\frac{1}{\alpha}}, &\alpha\in(1,2),\\
(\hat{D}_{\alpha}+D_{\alpha})n^{-1}\log n, &\alpha=1,\\
(\hat{D}_{\alpha}+D_{\alpha})n^{-1}, &\alpha\in(0,1).
\end{cases}
\end{align}
Now, let us bound the $\mathbf{I}$ and $\mathbf{II}.$

\vspace{2mm}
\noindent i) When $\alpha\in(1,2),$ one has by Lemma \ref{ml},
\begin{align*}
\mathbf{I}\leq\frac{4\|f''\|_{\infty}}{(2-\alpha)\sigma^{2}}(2A)^{\frac{2}{\alpha}}n^{\frac{\alpha-2}{\alpha}}+\frac{8\|f'\|_{\infty}}{(\alpha-1)\sigma^{\alpha}}\sup_{|x|\geq\sigma n^{\frac{1}{\alpha}}}|\epsilon(x)|+\frac{2\|f''\|_{\infty}}{\sigma^{2}}n^{\frac{\alpha-2}{\alpha}}\int_{-\sigma n^{\frac{1}{\alpha}}}^{\sigma n^{\frac{1}{\alpha}}}\frac{|\epsilon(x)|}{x^{\alpha-1}}\dif x,
\end{align*}
whereas
\begin{align*}
|\mathbf{II}|\leq\|f''\|_{\infty}\sigma^{-2}n^{-\frac{2}{\alpha}}\sum_{i=1}^{n}|\mathbb{E}[X_{i}]|\big(\mathbb{E}|X_{i}|+|\mathbb{E}[X_{i}]|\big)\leq\frac{2}{\sigma^{2}}\mathbb{E}|X_{1}||\mathbb{E}[X_{1}]| \|f''\|_{\infty}n^{-\frac{2-\alpha}{\alpha}}.
\end{align*}

\vspace{2mm}
\noindent ii) When $\alpha=1,$ $\beta=0$ and $\gamma\in(0,\infty),$ one has by Lemma \ref{ml1},
\begin{align*}
\mathbf{I}\leq\frac{12A^{2}\|f''\|_{\infty}}{\sigma}n^{-1}&+\frac{2(2\|f\|_{\infty}+\|f'\|_{\infty})}{\sigma}\sup_{|x|\geq \sigma n}|\epsilon(x)|\\
&+\frac{\|f''\|_{\infty}}{\sigma^{2}}n^{-1}\int_{-\sigma n}^{\sigma n}|\epsilon(x)|dx+\frac{\|f'\|_{\infty}}{\sigma}\int_{\sigma n}^{\infty}\frac{|\epsilon(x)|+|\epsilon(-x)|}{x}dx.
\end{align*}
whereas by (\ref{expectation}) and Lemma \ref{random truncation},
\begin{align*}
|\mathbf{II}|\leq&\frac{n^{-1}}{\sigma}\sum_{i=1}^{n}\Big\{\mathbb{E}\Big|\int_{0}^{1}\mathbb{E}[X_{i}{\bf 1}_{(0,\sigma n)}(|X_{i}|)]\Big[f'\Big(\frac{Z_{i}}{n}\Big)-f'\Big(\frac{Z_{i}}{n}+\frac{X_{i}}{\sigma n}\Big)\Big]\dif t\Big|\\
&+\mathbb{E}\Big|\int_{0}^{1}\mathbb{E}[X_{i}{\bf 1}_{(0,\sigma n)}(|X_{i}|)]\Big[f'\Big(\frac{Z_{i}}{n}+\frac{X_{i}}{\sigma n}\Big)-f'\Big(\frac{Z_{i}}{n}+\frac{X_{i}}{\sigma n}-\frac{n^{-1}}{\sigma}\mathbb{E}\big[X_{i}{\bf 1}_{(0,\sigma n)}(|X_{i}|)\big]t\Big)\Big]\dif t\Big|\Big\}\\
\leq&\frac{4}{\sigma^{2}}\big(\frac{K}{\gamma}+k+1\big)\Big[\|f''\|_{\infty}\big(\frac{K}{\gamma}+K+2\big)+(A+K)\big(2\|f'\|_{\infty}+\|f''\|_{\infty}\log(\sigma n)\big)\Big]n^{-1}.
\end{align*}

\vspace{2mm}
\noindent iii) When $\alpha\in(0,1),$ on the one hand, using Lemma \ref{ml2},

a.) When $\gamma\in(1-\alpha,\infty),$ we have
\begin{align*}
\mathbf{I}\leq\frac{2+\alpha}{(2-\alpha)\sigma^{2}}&(2A)^{\frac{2}{\alpha}}\|f''\|_{\infty}n^{\frac{\alpha-2}{\alpha}}+\frac{2(3\|f\|_{\infty}
+2\|f'\|_{\infty})}{\sigma^{\alpha}}\sup_{|x|\geq\sigma n^{\frac{1}{\alpha}}}|\epsilon(x)|\\&+\frac{2\|f''\|_{\infty}}{\sigma^{2}}n^{\frac{\alpha-2}{\alpha}}\int_{-\sigma n^{\frac{1}{\alpha}}}^{\sigma n^{\frac{1}{\alpha}}}\frac{|\epsilon(x)|}{|x|^{\alpha-1}}dx+\frac{\|f'\|_{\infty}}{\sigma}n^{\frac{\alpha-1}{\alpha}}\int_{\sigma n^{\frac{1}{\alpha}}}^{\infty}\frac{|\epsilon(x)|+|\epsilon(-x)|}{x^{\alpha}}dx.
\end{align*}

b.) When $\gamma\in[0,1-\alpha],$ we have
\begin{align*}
\mathbf{I}\leq\frac{2+\alpha}{(2-\alpha)\sigma^{2}}&(2A)^{\frac{2}{\alpha}}\|f''\|_{\infty}n^{\frac{\alpha-2}{\alpha}}+\frac{2\|f''\|_{\infty}}{\sigma^{2}}n^{\frac{\alpha-2}{\alpha}}\int_{-\sigma n^{\frac{1}{\alpha}}}^{\sigma n^{\frac{1}{\alpha}}}\frac{|\epsilon(x)|}{|x|^{\alpha-1}}dx\\
&+\big[(4A+6K+4)\|f\|_{\infty}+\frac{8-4\alpha}{1-\alpha}\|f'\|_{\infty}\big]\sigma^{-\alpha}\sup_{|x|\geq\sigma n^{\frac{1}{\alpha}}}|\epsilon(x)|^{\alpha}.
\end{align*}
On the other hand, by Lemma \ref{expectation} we have
\begin{align*}
\mathbb{E}\Big[X_{i}{\bf 1}_{(0,\sigma n^{\frac{1}{\alpha}})}(|X_{i}|)\Big]=\frac{2A\alpha\beta}{1-\alpha}\big(\frac{n^{-\frac{1}{\alpha}}}{\sigma}\big)^{\alpha-1}&+(1+\beta)\int_{0}^{\sigma n^{\frac{1}{\alpha}}}\frac{\epsilon(x)}{x^{\alpha}}\dif x-(1-\beta)\int_{0}^{\sigma n^{\frac{1}{\alpha}}}\frac{\epsilon(-x)}{x^{\alpha}}\dif x\\
&+\big(\frac{n^{-\frac{1}{\alpha}}}{\sigma}\big)^{\alpha-1}\big[(1-\beta)\epsilon(-\sigma n^{\frac{1}{\alpha}})-(1+\beta)\epsilon(\sigma n^{\frac{1}{\alpha}})\big],
\end{align*}
which follows that
\begin{align*}
\mathbf{II}\leq\|f'\|_{\infty}\Big(\frac{4}{\sigma^{\alpha}}\sup_{|x|\geq\sigma n^{\frac{1}{\alpha}}}|\epsilon(x)|+\frac{n^{\frac{\alpha-1}{\alpha}}}{\sigma}\Big|(1+\beta)\int_{0}^{\sigma n^{\frac{1}{\alpha}}}\frac{\epsilon(x)}{x^{\alpha}}\dif x-(1-\beta)\int_{0}^{\sigma n^{\frac{1}{\alpha}}}\frac{\epsilon(-x)}{x^{\alpha}}\dif x\Big|\Big).
\end{align*}

Combining all of above, we get the desired conclusion of Theorem \ref{main}.
\end{proof}

\begin{proof} [{\bf Proof of Theorem \ref{improve}}]
It suffices to bound the $\mathbf{I}$ in the proof of Theorem \ref{main}. By the same argument as (\ref{1result}), we have
\begin{align*}
\mathbf{I}\leq\sum_{i=1}^{n}\Big\{\mathbb{E}\Big|\int_{-\infty}^{\infty}\Big[f\big(\frac{Z_{i}}{n^{\frac{1}{\alpha}}}+\frac{n^{-\frac{1}{\alpha}}}{\sigma}x\big)-\frac{n^{-\frac{1}{\alpha}}}{\sigma}x{\bf 1}_{(-1,1)}\big(\frac{n^{-\frac{1}{\alpha}}}{\sigma}x)f'\big(\frac{Z_{i}}{n^{\frac{1}{\alpha}}}\big)\Big]d\big[F_{X}(x)-F_{\tilde{X}}(x)\big]\Big|+\mathcal{R}_{i}\Big\},
\end{align*}
where
\begin{align*}
\mathcal{R}_{i}&=2A\alpha\mathbb{E}\Big|\int_{|x|<(2A)^{\frac{1}{\alpha}}}\frac{f\big(\frac{Z_{i}}{n^{\frac{1}{\alpha}}}+\frac{n^{-\frac{1}{\alpha}}}{\sigma}x\big)-f\big(\frac{Z_{i}}{n^{\frac{1}{\alpha}}}\big)
-\frac{n^{-\frac{1}{\alpha}}}{\sigma}xf'\big(\frac{Z_{i}}{n^{\frac{1}{\alpha}}}\big)}{|x|^{1+\alpha}}dx\Big|\\
&\leq\frac{4A\alpha}{(2-\alpha)\sigma^{2}}(2A)^{\frac{2-\alpha}{\alpha}}\|f''\|_{\infty}n^{-\frac{2}{\alpha}}.
\end{align*}
For the first term, we have by the same argument as (\ref{e:IBP})
\begin{align*}
&\mathbb{E}\Big|\int_{-\infty}^{\infty}\Big[f\big(\frac{Z_{i}}{n^{\frac{1}{\alpha}}}+\frac{n^{-\frac{1}{\alpha}}}{\sigma}x\big)-\frac{n^{-\frac{1}{\alpha}}}{\sigma}x{\bf 1}_{(-1,1)}\big(\frac{n^{-\frac{1}{\alpha}}}{\sigma}x)f'\big(\frac{Z_{i}}{n^{\frac{1}{\alpha}}}\big)\Big]d\big[F_{X}(x)-F_{\tilde{X}}(x)\big]\Big|\\
\leq&\frac{2(2A)^{\frac{2}{\alpha}}}{\sigma^{2}}\|f''\|_{\infty}n^{-\frac{2}{\alpha}}+\frac{4\|f\|_{\infty}}{\sigma^{\alpha}}n^{-1}\sup_{|x|\geq\sigma n^{\frac{1}{\alpha}}}|\epsilon(x)|+\frac{2\|f''\|_{\infty}}{\sigma^{2}}n^{-\frac{2}{\alpha}}\int_{-\sigma n^{\frac{1}{\alpha}}}^{\sigma n^{\frac{1}{\alpha}}}\frac{|\epsilon(x)|}{|x|^{\alpha-1}}dx,
\end{align*}
the desired conclusion follows.
\end{proof}

\section{A more difficult example: Proof of (\ref{ex}) and (\ref{ex1})}\label{difficult}
In this section, we prove the estimates (\ref{ex}) and (1.9).
Consider independent copies $X_1,\ldots,X_n$ of a random variable $X$ with density
\begin{align*}
P(X>x)=K_{0}\frac{\big(\log|x|\big)^{\delta}}{|x|^{\alpha}},\quad x\geq e,\qquad P(X<x)=K_{0}\frac{\big(\log|x|\big)^{\delta}}{|x|^{\alpha}},\quad x\leq-e,
\end{align*}
where $K_{0}>0,$ $\alpha\in(0,2)$ and $\delta\in\mathbb{R}.$ The corresponding density function is
\begin{align*}
p_{X}(x)=\frac{K_{0}\big[\alpha(\log|x|)^{\delta}-\delta(\log|x|)^{\delta-1}\big]}{|x|^{\alpha+1}}{\bf 1}_{[e,\infty)}(|x|).
\end{align*}
Recall $\gamma_{n}=\inf\{x>0:\mathbb{P}(|X|>x)\leq\frac{1}{n}\}$ can be determined by
\begin{align}\label{log}
\frac{n}{\gamma_{n}^{\alpha}}=\frac{1}{2K_{0}(\log\gamma_{n})^{\delta}}
\end{align}
and it is easy to see $C_{\alpha,\delta}n^{\frac{1}{\alpha}}\leq\gamma_{n}\leq C_{\alpha,\delta}n^{\frac{1}{\alpha}}(\log\gamma_{n})^{\delta}.$ Now, we set $\widetilde{X}_i=\frac{n^\frac1\alpha}{\tilde{\sigma}\gamma_n}X_{i}$ and $\tilde{Z}_{i}=\tilde{Y}_{1}+\cdots+\tilde{Y}_{i-1}+\tilde{X}_{i+1}+\cdots+\tilde{X}_{n}$  for $1\leq i\leq n,$ where $\tilde{Y}_{1},\cdots,\tilde{Y}_{n}$ are independent copies of $\tilde{Y}.$

By the same argument as (\ref{aL}), for any $f\in\mathcal{C}_{b}^{3}(\mathbb{R}),$ we have
\begin{align*}
\mathbb{E}[f(\tilde{S}_{n})]-\mathbb{E}[f(\tilde{Y})]=\sum_{i=1}^{n}\mathbb{E}\Big[f\Big(\frac{\tilde{X}_{i}+\tilde{Z}_{i}}{n^{\frac{1}{\alpha}}}\Big)-f\Big(\frac{\tilde{Z}_{i}}{n^{\frac{1}{\alpha}}}\Big)\Big]-
\sum_{i=1}^{n}\mathbb{E}\Big[f\Big(\frac{\tilde{Y}_{i}+\tilde{Z}_{i}}{n^{\frac{1}{\alpha}}}\Big)-f\Big(\frac{\tilde{Z}_{i}}{n^{\frac{1}{\alpha}}}\Big)\Big].
\end{align*}
Using (\ref{log}), we have
\begin{align*}
&\quad\mathbb{E}\Big[f\Big(\frac{\tilde{X}_{i}+\tilde{Z}_{i}}{n^{\frac{1}{\alpha}}}\Big)-f\Big(\frac{\tilde{Z}_{i}}{n^{\frac{1}{\alpha}}}\Big)\Big]\\
&=\frac{d_{\alpha}}{2n}\mathbb{E}\Big\{\int_{\mathbb{R}}\Big[f\Big(\frac{\tilde{Z}_{i}}{n^{\frac{1}{\alpha}}}+u\Big)
-f\Big(\frac{\tilde{Z}_{i}}{n^{\frac{1}{\alpha}}}\Big)\Big]\frac{\alpha\big((\log(\tilde{\sigma}\gamma_{n}|u|)\big)^{\delta}-\delta\big((\log(\tilde{\sigma}\gamma_{n}|u|)\big)^{\delta-1}}{\alpha(\log\gamma_{n})^{\delta}|u|^{\alpha+1}}{\bf 1}_{[e,\infty)}(\tilde{\sigma}\gamma_{n}|u|)\dif u\Big\}\\
&=\frac{d_{\alpha}}{2n}\mathbb{E}\Big\{\int_{\mathbb{R}}\Big[f\Big(\frac{\tilde{Z}_{i}}{n^{\frac{1}{\alpha}}}+u\Big)
-f\Big(\frac{\tilde{Z}_{i}}{n^{\frac{1}{\alpha}}}\Big)\Big]\frac{1}{|u|^{\alpha+1}}{\bf 1}_{[e,\infty)}(\tilde{\sigma}\gamma_{n}|u|)\dif u\Big\}\\
&\quad+\frac{d_{\alpha}}{2n}\mathbb{E}\Big\{\int_{\mathbb{R}}\Big[f\Big(\frac{\tilde{Z}_{i}}{n^{\frac{1}{\alpha}}}+u\Big)
-f\Big(\frac{\tilde{Z}_{i}}{n^{\frac{1}{\alpha}}}\Big)\Big]\frac{\big(\log(\tilde{\sigma}\gamma_{n}|u|)\big)^{\delta}-\big(\log(\gamma_{n})\big)^{\delta}}{(\log\gamma_{n})^{\delta}|u|^{\alpha+1}}{\bf 1}_{[e,\infty)}(\tilde{\sigma}\gamma_{n}|u|)\dif u\Big\}\\
&\quad-\frac{d_{\alpha}}{2n}\mathbb{E}\Big\{\int_{\mathbb{R}}\Big[f\Big(\frac{\tilde{Z}_{i}}{n^{\frac{1}{\alpha}}}+u\Big)
-f\Big(\frac{\tilde{Z}_{i}}{n^{\frac{1}{\alpha}}}\Big)\Big]\frac{\delta\big((\log(\tilde{\sigma}\gamma_{n}|u|)\big)^{\delta-1}}{\alpha(\log\gamma_{n})^{\delta}|u|^{\alpha+1}}{\bf 1}_{[e,\infty)}(\tilde{\sigma}\gamma_{n}|u|)\dif u\Big\}.
\end{align*}
On the one hand, we have
\begin{align*}
\frac{1}{n}\mathbb{E}\Big[\mathcal{L}^{\alpha,0}f\Big(\frac{\tilde{Z}_{i}}{n^{\frac{1}{\alpha}}}\Big)\Big]=\frac{d_{\alpha}}{2n}\mathbb{E}\Big\{\int_{\mathbb{R}}\Big[f\Big(\frac{\tilde{Z}_{i}}{n^{\frac{1}{\alpha}}}+u\Big)
-f\Big(\frac{\tilde{Z}_{i}}{n^{\frac{1}{\alpha}}}\Big)\Big]\frac{1}{|u|^{\alpha+1}}\dif u\Big\},
\end{align*}
\begin{align*}
\mathbb{E}\Big[f\Big(\frac{\tilde{Y}_{i}+\tilde{Z}_{i}}{n^{\frac{1}{\alpha}}}\Big)-f\Big(\frac{\tilde{Z}_{i}}{n^{\frac{1}{\alpha}}}\Big)\Big]=\frac{1}{n}\int_{0}^{1}\mathbb{E}\Big[\mathcal{L}^{\alpha,0}
f\Big(\frac{\tilde{Z}_{i}+\hat{Y}_{s}}{n^{\frac{1}{\alpha}}}\Big)\Big]\dif s.
\end{align*}
As a result,
\begin{align*}
|\mathbb{E}[f(\tilde{S}_{n})]-\mathbb{E}[f(\tilde{Y})]|\leq\frac{1}{n}\sum_{i=1}^{n}\int_{0}^{1}\mathbb{E}\Big|\mathcal{L}^{\alpha,0}f\Big(\frac{\tilde{Z}_{i}}{n^{\frac{1}{\alpha}}}\Big)
-\mathcal{L}^{\alpha,0}f\Big(\frac{\tilde{Z}_{i}+\hat{Y}_{s}}{n^{\frac{1}{\alpha}}}\Big)\Big|\dif s+\mathcal{I}+\mathcal{II}+\mathcal{III},
\end{align*}
where
\begin{align*}
\mathcal{I}:=\frac{d_{\alpha}}{2n}\sum_{i=1}^{n}\mathbb{E}\Big|\int_{\mathbb{R}}\Big[f\Big(\frac{\tilde{Z}_{i}}{n^{\frac{1}{\alpha}}}+u\Big)
-f\Big(\frac{\tilde{Z}_{i}}{n^{\frac{1}{\alpha}}}\Big)\Big]\frac{1}{|u|^{\alpha+1}}{\bf 1}_{[0,e)}(\tilde{\sigma}\gamma_{n}|u|)\dif u\Big|,
\end{align*}
\begin{align*}
\mathcal{II}:=\frac{d_{\alpha}}{2n}\sum_{i=1}^{n}\mathbb{E}\Big|\int_{\mathbb{R}}\Big[f\Big(\frac{\tilde{Z}_{i}}{n^{\frac{1}{\alpha}}}+u\Big)
-f\Big(\frac{\tilde{Z}_{i}}{n^{\frac{1}{\alpha}}}\Big)\Big]\frac{\big(\log(\tilde{\sigma}\gamma_{n}|u|)\big)^{\delta}-\big(\log(\gamma_{n})\big)^{\delta}}{(\log\gamma_{n})^{\delta}|u|^{\alpha+1}}{\bf 1}_{[e,\infty)}(\tilde{\sigma}\gamma_{n}|u|)\dif u\Big|
\end{align*}
and
\begin{align*}
\mathcal{III}:=\frac{d_{\alpha}}{2n}\sum_{i=1}^{n}\mathbb{E}\Big|\int_{\mathbb{R}}\Big[f\Big(\frac{\tilde{Z}_{i}}{n^{\frac{1}{\alpha}}}+u\Big)
-f\Big(\frac{\tilde{Z}_{i}}{n^{\frac{1}{\alpha}}}\Big)\Big]\frac{\delta\big((\log(\tilde{\sigma}\gamma_{n}|u|)\big)^{\delta-1}}{\alpha(\log\gamma_{n})^{\delta}|u|^{\alpha+1}}{\bf 1}_{[e,\infty)}(\tilde{\sigma}\gamma_{n}|u|)\dif u\Big|.
\end{align*}

By (\ref{log1}) with $\beta=0$ and $Z_{i}$ replaced by $\tilde{Z}_{i},$ we know
\begin{align*}
\frac{1}{n}\sum_{i=1}^{n}\int_{0}^{1}\mathbb{E}\Big|\mathcal{L}^{\alpha,0}f\Big(\frac{\tilde{Z}_{i}}{n^{\frac{1}{\alpha}}}\Big)
-\mathcal{L}^{\alpha,0}f\Big(\frac{\tilde{Z}_{i}+\hat{Y}_{s}}{n^{\frac{1}{\alpha}}}\Big)\Big|\dif s
\leq C_{\alpha}\begin{cases}
D_{\alpha}n^{-\frac{1}{\alpha}}, &\alpha\in(1,2),\\
(\hat{D}_{\alpha}+D_{\alpha})n^{-1}\log n, &\alpha=1,\\
(\hat{D}_{\alpha}+D_{\alpha})n^{-1}, &\alpha\in(0,1).
\end{cases}
\end{align*}
On the other hand,
\begin{align*}
\mathcal{I}&=\frac{d_{\alpha}}{2n}\sum_{i=1}^{n}\mathbb{E}\Big|\int_{-(\tilde{\sigma}\gamma_{n})^{-1}}^{(\tilde{\sigma}\gamma_{n})^{-1}}\Big[f\Big(\frac{\tilde{Z}_{i}}{n^{\frac{1}{\alpha}}}+u\Big)
-f\Big(\frac{\tilde{Z}_{i}}{n^{\frac{1}{\alpha}}}\Big)-uf'\Big(\frac{\tilde{Z}_{i}}{n^{\frac{1}{\alpha}}}\Big)\Big]\frac{1}{|u|^{\alpha+1}}\dif u\Big|\\
&=\frac{d_{\alpha}}{2n}\sum_{i=1}^{n}\mathbb{E}\Big|\int_{-(\tilde{\sigma}\gamma_{n})^{-1}}^{(\tilde{\sigma}\gamma_{n})^{-1}}\int_{0}^{1}
u\Big[f'\Big(\frac{\tilde{Z}_{i}}{n^{\frac{1}{\alpha}}}
+ut\Big)-f'\Big(\frac{\tilde{Z}_{i}}{n^{\frac{1}{\alpha}}}\Big)\Big]\frac{1}{|u|^{\alpha+1}}\dif t\dif u\Big|\\
&\leq C_{\alpha}\|f''\|_{\infty}\int_{0}^{(\tilde{\sigma}\gamma_{n})^{-1}}\frac{\dif u}{u^{\alpha-1}}=O\big(\gamma_{n}^{\alpha-2}\big).
\end{align*}
For $\mathcal{II}$, we choose a number $\tilde{\sigma}^{-1}<N<\gamma_{n}$ and divide $\mathcal{II}$ into the following two terms:
\begin{align*}
\mathcal{II}&\leq \frac{d_{\alpha}}{2n}\sum_{i=1}^{n}\mathbb{E}\Big|\int_{\mathbb{R}}\Big[f\Big(\frac{\tilde{Z}_{i}}{n^{\frac{1}{\alpha}}}+u\Big)
-f\Big(\frac{\tilde{Z}_{i}}{n^{\frac{1}{\alpha}}}\Big)\Big]\frac{\big(\log(\tilde{\sigma}\gamma_{n}|u|)\big)^{\delta}-\big(\log(\gamma_{n})\big)^{\delta}}{(\log\gamma_{n})^{\delta}|u|^{\alpha+1}}{\bf 1}_{[(\tilde{\sigma}\gamma_{n})^{-1}e,N)}(|u|)\dif u\Big|\\
&\quad+\frac{d_{\alpha}}{2n}\sum_{i=1}^{n}\mathbb{E}\Big|\int_{\mathbb{R}}\Big[f\Big(\frac{\tilde{Z}_{i}}{n^{\frac{1}{\alpha}}}+u\Big)
-f\Big(\frac{\tilde{Z}_{i}}{n^{\frac{1}{\alpha}}}\Big)\Big]\frac{\big(\log(\tilde{\sigma}\gamma_{n}|u|)\big)^{\delta}-\big(\log(\gamma_{n})\big)^{\delta}}{(\log\gamma_{n})^{\delta}|u|^{\alpha+1}}{\bf 1}_{[N,\infty)}(|u|)\dif u\Big|\\
&:=\mathcal{II}_{1}+\mathcal{II}_{2}.
\end{align*}
One has
\begin{align*}
\mathcal{II}_{1}&=\frac{d_{\alpha}}{2n}\sum_{i=1}^{n}\mathbb{E}\Big|\int_{\mathbb{R}}\Big[f\Big(\frac{\tilde{Z}_{i}}{n^{\frac{1}{\alpha}}}+u\Big)
-f\Big(\frac{\tilde{Z}_{i}}{n^{\frac{1}{\alpha}}}\Big)-uf'\Big(\frac{\tilde{Z}_{i}}{n^{\frac{1}{\alpha}}}\Big)\Big]\\
&\qquad\qquad\qquad\qquad\qquad\qquad\qquad\qquad\cdot\frac{\big(\log(\tilde{\sigma}\gamma_{n}|u|)\big)^{\delta}-\big(\log(\gamma_{n})\big)^{\delta}}{(\log\gamma_{n})^{\delta}|u|^{\alpha+1}}{\bf 1}_{[(\tilde{\sigma}\gamma_{n})^{-1}e,N)}(|u|)\dif u\Big|\\
&\leq C_{\alpha}n^{-1}\!\sum_{i=1}^{n}\mathbb{E}\int_{\mathbb{R}}\Big|f\Big(\frac{\tilde{Z}_{i}}{n^{\frac{1}{\alpha}}}+u\Big)
-f\Big(\frac{\tilde{Z}_{i}}{n^{\frac{1}{\alpha}}}\Big)-uf'\Big(\frac{\tilde{Z}_{i}}{n^{\frac{1}{\alpha}}}\Big)\Big|\frac{\big|\big(1+\frac{\log\tilde{\sigma}|u|}{\log\gamma_{n}}\big)^{\delta}-1\big|}{|u|^{\alpha+1}}{\bf 1}_{[(\tilde{\sigma}\gamma_{n})^{-1}e,N)}(|u|)\dif u\\
&\leq C_{\alpha}
\begin{cases}
\|f'\|_{\infty}\int_{\frac{e}{\tilde{\sigma}\gamma_{n}}}^{N}
\frac{\big|\log u\big|}{\log\gamma_{n}u^{\alpha}}\dif u, &\alpha\in(0,1)\\
\|f''\|_{\infty}\int_{\frac{e}{\tilde{\sigma}\gamma_{n}}}^{N}\frac{\big|\log u\big|}{\log\gamma_{n}u^{\alpha-1}}\dif u, &\alpha\in[1,2)
\end{cases}\\
&\leq C_{\alpha}
\begin{cases}
\|f'\|_{\infty}\frac{N^{1-\alpha}\log N}{\log\gamma_{n}}, &\alpha\in(0,1)\\
\|f''\|_{\infty}\frac{N^{2-\alpha}\log N}{\log\gamma_{n}}, &\alpha\in[1,2),
\end{cases}
\end{align*}
the last second inequality is by the fact $\gamma_{n}^{-1}e\leq\tilde{\sigma}|u|<\tilde{\sigma}N\leq\tilde{\sigma}\gamma_{n}$ and $\big||1+x|^{\delta}-1\big|\leq C_{\delta}|x|$ for any $|x|<1.$
For the second term, when $\delta\in(-\infty,0],$ we have
\begin{align*}
\mathcal{II}_{2}&=\frac{d_{\alpha}}{2n}\sum_{i=1}^{n}\mathbb{E}\Big|\int_{\mathbb{R}}\Big[f\Big(\frac{\tilde{Z}_{i}}{n^{\frac{1}{\alpha}}}+u\Big)
-f\Big(\frac{\tilde{Z}_{i}}{n^{\frac{1}{\alpha}}}\Big)\Big]\frac{\big(\log(\tilde{\sigma}\gamma_{n}|u|)\big)^{\delta}-\big(\log(\gamma_{n})\big)^{\delta}}{(\log\gamma_{n})^{\delta}|u|^{\alpha+1}}{\bf 1}_{[N,\infty)}(|u|)\dif u\Big|\\
&\leq C_{\alpha}n^{-1}\sum_{i=1}^{n}\mathbb{E}\int_{\mathbb{R}}\Big|f\Big(\frac{\tilde{Z}_{i}}{n^{\frac{1}{\alpha}}}+u\Big)
-f\Big(\frac{\tilde{Z}_{i}}{n^{\frac{1}{\alpha}}}\Big)\Big|\frac{(\log\gamma_{n})^{\delta}}{(\log\gamma_{n})^{\delta}|u|^{\alpha+1}}{\bf 1}_{[N,\infty)}(|u|)\dif u\\
&\leq C_{\alpha}
\begin{cases}
\|f\|_{\infty}N^{-\alpha}, &\alpha\in(0,1],\\
\|f'\|_{\infty}N^{1-\alpha}, &\alpha\in(1,2).
\end{cases}
\end{align*}
In addition, according to the $C_{r}$-inequality: for any $a,b\in\mathbb{R}$ and $r>0,$
$$
|a+b|^{r}\leq\max\{1,2^{r-1}\}\big(|a|^{r}+|b|^{r}\big),
$$
when $\delta\in(0,1],$ we have
\begin{align}\label{need}
\mathcal{II}_{2}&\leq C_{\alpha}n^{-1}\sum_{i=1}^{n}\mathbb{E}\int_{\mathbb{R}}\Big|f\Big(\frac{\tilde{Z}_{i}}{n^{\frac{1}{\alpha}}}+u\Big)
-f\Big(\frac{\tilde{Z}_{i}}{n^{\frac{1}{\alpha}}}\Big)\Big|\frac{\big(\log(\tilde{\sigma}|u|)\big)^{\delta}}{(\log\gamma_{n})^{\delta}|u|^{\alpha+1}}{\bf 1}_{[N,\infty)}(|u|)\dif u\nonumber\\
&\leq C_{\alpha,\delta}
\begin{cases}
\|f\|_{\infty}\frac{(\log N)^{\delta}}{(\log\gamma_{n})^{\delta}}N^{-\alpha}, &\alpha\in(0,1],\\
\|f'\|_{\infty}\frac{(\log N)^{\delta}}{(\log\gamma_{n})^{\delta}}N^{1-\alpha}, &\alpha\in(1,2).
\end{cases}\\
&\leq C_{\alpha,\delta}
\begin{cases}
\|f\|_{\infty}N^{-\alpha}, &\alpha\in(0,1],\\
\|f'\|_{\infty}N^{1-\alpha}, &\alpha\in(1,2).
\end{cases}\nonumber
\end{align}
When $\delta\in(1,\infty),$ we have
\begin{align*}
\mathcal{II}_{2}&\leq C_{\alpha}n^{-1}\sum_{i=1}^{n}\mathbb{E}\int_{\mathbb{R}}\Big|f\Big(\frac{\tilde{Z}_{i}}{n^{\frac{1}{\alpha}}}+u\Big)
-f\Big(\frac{\tilde{Z}_{i}}{n^{\frac{1}{\alpha}}}\Big)\Big|\frac{\big(\log(\tilde{\sigma}|u|)\big)^{\delta}+(\log\gamma_{n})^{\delta}}{(\log\gamma_{n})^{\delta}|u|^{\alpha+1}}{\bf 1}_{[N,\infty)}(|u|)\dif u\\
&\leq C_{\alpha,\delta}
\begin{cases}
\|f\|_{\infty}N^{-\alpha}, &\alpha\in(0,1],\\
\|f'\|_{\infty}N^{1-\alpha}, &\alpha\in(1,2).
\end{cases}
\end{align*}
Hence, we let
\begin{align*}
\frac{N}{\log\gamma_{n}}=
\begin{cases}
N^{-\alpha}, &\alpha\in(0,1],\\
N^{1-\alpha}, &\alpha\in(1,2),
\end{cases}
\end{align*}
this implies
\begin{align*}
N=
\begin{cases}
(\log\gamma_{n})^{\frac{1}{\alpha+1}}, &\alpha\in(0,1],\\
(\log\gamma_{n})^{\frac{1}{\alpha}}, &\alpha\in(1,2),
\end{cases}
\end{align*}
which follow that
\begin{align*}
\mathcal{II}=
\begin{cases}
O\big((\log\gamma_{n})^{-\frac{\alpha}{\alpha+1}}\big), &\alpha\in(0,1),\\
O\big((\log\gamma_{n})^{\theta-\frac{1}{2}}\big), &\alpha=1,\\
O\big((\log\gamma_{n})^{-\frac{\alpha-1}{\alpha}}\big), &\alpha\in(1,2),
\end{cases}
=
\begin{cases}
O\big((\log n)^{-1+\frac{1}{\alpha+1}}\big), &\alpha\in(0,1),\\
O\big((\log n)^{\theta-\frac{1}{2}}\big), &\alpha=1,\\
O\big((\log n)^{-1+\frac{1}{\alpha}}\big), &\alpha\in(1,2),
\end{cases}
\end{align*}
for some very small $\theta>0.$\\
For $\mathcal{III},$ since
\begin{align*}
\mathcal{III}\leq&\frac{d_{\alpha}}{2n}\sum_{i=1}^{n}\mathbb{E}\Big|\int_{\mathbb{R}}\Big[f\Big(\frac{\tilde{Z}_{i}}{n^{\frac{1}{\alpha}}}+u\Big)
-f\Big(\frac{\tilde{Z}_{i}}{n^{\frac{1}{\alpha}}}\Big)\Big]\frac{\delta\big((\log(\tilde{\sigma}\gamma_{n}|u|)\big)^{\delta-1}-\delta(\log\gamma_{n})^{\delta-1}}{\alpha(\log\gamma_{n})^{\delta}|u|^{\alpha+1}}{\bf 1}_{[e,\infty)}(\tilde{\sigma}\gamma_{n}|u|)\dif u\Big|\\
&+\frac{d_{\alpha}}{2n}\sum_{i=1}^{n}\mathbb{E}\Big|\int_{\mathbb{R}}\Big[f\Big(\frac{\tilde{Z}_{i}}{n^{\frac{1}{\alpha}}}+u\Big)
-f\Big(\frac{\tilde{Z}_{i}}{n^{\frac{1}{\alpha}}}\Big)\Big]\frac{\delta(\log\gamma_{n})^{\delta-1}}{\alpha(\log\gamma_{n})^{\delta}|u|^{\alpha+1}}{\bf 1}_{[e,\infty)}(\tilde{\sigma}\gamma_{n}|u|)\dif u\Big|,
\end{align*}
and
\begin{align*}
\frac{d_{\alpha}}{2n}\sum_{i=1}^{n}\mathbb{E}\Big|\int_{\mathbb{R}}\Big[f\Big(\frac{\tilde{Z}_{i}}{n^{\frac{1}{\alpha}}}+u\Big)
-f\Big(\frac{\tilde{Z}_{i}}{n^{\frac{1}{\alpha}}}\Big)\Big]\frac{\delta(\log\gamma_{n})^{\delta-1}}{\alpha(\log\gamma_{n})^{\delta}|u|^{\alpha+1}}{\bf 1}_{[e,\infty)}(\tilde{\sigma}\gamma_{n}|u|)\dif u\Big|=O\big((\log n)^{-1}\big),
\end{align*}
we have by the same argument as the proof of $\mathcal{II},$
\begin{align*}
\mathcal{III}=O\big((\log n)^{-1}\big).
\end{align*}
Therefore, putting everything together, we get that
\begin{align*}
\big|\mathbb{E}[f(\tilde{S}_{n})]-\mathbb{E}[f(\tilde{Y})]\big|=
\begin{cases}
O\big((\log n)^{-1+\frac{1}{\alpha+1}}\big), &\alpha\in(0,1),\\
O\big((\log n)^{\theta-\frac{1}{2}}\big), &\alpha=1,\\
O\big((\log n)^{-1+\frac{1}{\alpha}}\big), &\alpha\in(1,2),
\end{cases}
\end{align*}
for some very small $\theta>0.$

In particular, when $\delta\in(0,1],$ we have by (\ref{need})
\begin{align*}
\mathcal{II}_{2}\leq C_{\alpha,\delta}
\begin{cases}
\|f\|_{\infty}\frac{(\log N)^{\delta}}{(\log\gamma_{n})^{\delta}}N^{-\alpha}, &\alpha\in(0,1],\\
\|f'\|_{\infty}\frac{(\log N)^{\delta}}{(\log\gamma_{n})^{\delta}}N^{1-\alpha}, &\alpha\in(1,2).
\end{cases}
\end{align*}
Taking $N=\tilde{\sigma}^{-1}+1,$ we have
\begin{align*}
\mathcal{II}=O\big((\log n)^{-\delta}\big).
\end{align*}
Putting everything together, we also have
\begin{align}\label{need2}
\big|\mathbb{E}[f(\tilde{S}_{n})]-\mathbb{E}[f(\tilde{Y})]\big|=O\big((\log n)^{-\delta}\big).
\end{align}
Combining (\ref{need2}) and (\ref{ex}), we immediately obtain (\ref{ex1}).
\qed

\br
For the sake of convenience, we let
\begin{align*}
\frac{N}{\log\gamma_{n}}=
\begin{cases}
N^{-\alpha}, &\alpha\in(0,1],\\
N^{1-\alpha}, &\alpha\in(1,2),
\end{cases}
\end{align*}
in the above proof. In fact. if we let $N\log N=\log\gamma_{n},$ then the rate will be better.
\er

\section{Appendix. Proof of Corollary \ref{c:SCLT}}\label{appendix}
\begin{proof}
We only give the proof of the case $\alpha\in[1,2)$ and the case $\alpha\in(0,1)$ is similar.

Consider and fix a function $f_{0}\in\mathcal{H}_{3}$ and $f_{0}:\mathbb{R}\rightarrow[0,1]$ such that $f_{0}(s)=1$ for $s\leq0$ and $f_{0}(s)=0$ for $s\geq1.$ Fix any $x\in\mathbb{R}$ and define $f(s)=f_{0}\big(\rho(s-x)\big)$ for some $\rho>1.$ For this function $f,$ we have
$$
\|f\|_{\infty}\leq C,\quad \|f'\|_{\infty}\leq C\rho,\quad \|f''\|_{\infty}\leq C\rho^{2}\quad and \quad\|f'''\|_{\infty}\leq C\rho^{3}
$$
for some constant $C>0.$ Then for any $y,z\in\mathbb{R},$ we have
\begin{align*}
\mathbb{P}(S_{n}\leq x)&\leq\mathbb{E}\big[f_{0}\big(\rho(S_{n}-x)\big)\big]=\mathbb{E}\big[f\big(S_{n})\big)\big]\\
&\leq\mathbb{E}\big[f(Y)\big]+4C\rho^{3}d_{\mcl W_3}(S_n,Y)\\
&\leq\mathbb{P}(Y\leq x+\frac{1}{\rho})+4C\rho^{3}d_{\mcl W_3}(S_n,Y),
\end{align*}
where the last second inequality thanks to Theorem \ref{main}. Furthermore, we have
\begin{align*}
\mathbb{P}(Y\leq x+\frac{1}{\rho})-\mathbb{P}(Y\leq x)\leq C_{\alpha}\int_{x}^{x+\frac{1}{\rho}}\frac{1}{(1+|y|)^{1+\alpha}}\dif y\leq C_{\alpha}\frac{1}{\rho}
\end{align*}
for some positive constant $C_{\alpha}$ depends on $\alpha.$ These imply
\begin{align*}
\mathbb{P}(S_{n}\leq x)-\mathbb{P}(Y\leq x)\leq4C\rho^{3}d_{\mcl W_3}(S_n,Y)+C_{\alpha}\frac{1}{\rho}.
\end{align*}
Hence, we take $\rho=\big[d_{\mcl W_3}(S_n,Y)\big]^{-\frac{1}{4}},$
\begin{align*}
\mathbb{P}(S_{n}\leq x)-\mathbb{P}(Y\leq x)=O\Big(\big(d_{\mcl W_3}(S_n,Y)\big)^{\frac{1}{4}}\Big).
\end{align*}
This gives one half of the claim. The other half follows similarly.
\end{proof}
\bigskip
{\bf Acknowledgements}: We would like to gratefully thank Persi Diaconis, Elton Hsu and Renming Song for very helpful discussions and suggestions. Special thanks are due to the two anonymous referees for carefully reviewing the paper and giving many useful suggestions. This work was partially supported by the grant 346300 for IMPAN from the Simons Foundation and the matching 2015-2019 Polish MNiSW fund. This research is also partly supported by the following grants:  Macao S.A.R. (FDCT 038/2017/A1, FDCT 030/2016/A1, FDCT 025/2016/A1), NNSFC 11571390, University of Macau MYRG (2016-00025-FST, 2018-00133-FST).

\end{document}